\documentclass[journal]{IEEEtran}
\newif\ifpdf
\ifx\pdfoutput\undefined
	\pdffalse
\else
	\pdfoutput=1
	\pdftrue
\fi
\ifpdf
	\usepackage{graphicx,graphics,latexsym,verbatim,epsfig}
	\usepackage{graphicx}
	\DeclareGraphicsExtensions{.pdf,.jpg,.png}
\else
	\usepackage{graphicx}
	\DeclareGraphicsExtensions{.eps}
\fi
\let\ifpdf\relax

\makeatletter
\def\maxwidth{
  \ifdim\Gin@nat@width>\linewidth
    \linewidth
  \else
    \Gin@nat@width
  \fi
}
\newdimen\mymathindent
\newenvironment{bulletequation}%
    {\@beginparpenalty\predisplaypenalty
     \@endparpenalty\postdisplaypenalty
     \refstepcounter{equation}%
     \trivlist \item[]\leavevmode
       \hb@xt@\linewidth\bgroup $\m@th
         \displaystyle
         \hskip\mymathindent}%
        {$\hfil 
         \displaywidth\linewidth\hbox{\@eqnnum}%
       \egroup
     \endtrivlist}
\makeatother

\usepackage{placeins}
\usepackage[cmex10]{amsmath}
\usepackage{amsfonts}
\usepackage{amsthm}
\usepackage{amssymb}
\usepackage{mathtools}
\usepackage{multirow}
\usepackage{enumerate}
\usepackage[caption=false,font=footnotesize]{subfig}
\usepackage{xr}
\newtheorem{theorem}{Theorem}

\newtheorem{lemma}{Lemma}

\newtheorem{proposition}{Proposition}

\newcommand{\refe}[1]{(\ref{#1})}

\newcommand{\Mzero}{\ensuremath \left.0\mkern-6.5mu \raisebox{0.6pt}{\text{\small /}}\right.}

\newcommand{\Floor}[1]{\left\lfloor #1 \right\rfloor}

\begin{document}
\title{A New Notion of Effective Resistance for Directed Graphs---Part II: Computing Resistances}
\author{George~Forrest~Young,~\IEEEmembership{Student~Member,~IEEE,
} Luca~Scardovi,~\IEEEmembership{Member,~IEEE,} and~Naomi~Ehrich~Leonard,~\IEEEmembership{Fellow,~IEEE}%
\thanks{This research was supported in part by AFOSR grant FA9550-07-1-0-0528, ONR grant N00014-09-1-1074, ARO grant W911NG-11-1-0385 and the Natural Sciences and Engineering Research Council (NSERC) of Canada.}%
\thanks{G. F. Young and N. E. Leonard are with the Department of Mechanical and Aerospace Engineering, Princeton University, Princeton, NJ 08544, USA; e-mail: \texttt{gfyoung@princeton.edu}; \texttt{naomi@princeton.edu}.}%
\thanks{L. Scardovi is with the Department of Electrical and Computer Engineering, University of Toronto, Toronto, ON, M5S 3G4, Canada; e-mail: \texttt{scardovi@scg.utoronto.ca}.}}

\date{\today}
\maketitle

\begin{abstract}
In Part I of this work we defined a generalization of the concept of effective resistance to directed graphs, and we explored some of the properties of this new definition. Here, we use the theory developed in Part I to compute effective resistances in some prototypical directed graphs. This exploration highlights cases where our notion of effective resistance for directed graphs behaves analogously to our experience from undirected graphs, as well as cases where it behaves in unexpected ways.
\end{abstract}

\begin{IEEEkeywords}
\noindent Graph theory, networks, networked control systems, directed graphs, effective resistance
\end{IEEEkeywords}

\section{Introduction}\label{sec:intro}
In the companion paper to this work, \cite{Young13I}, we presented a generalization of the concept of effective resistance to directed graphs. This extension was constructed algebraically to preserve the relationships for directed graphs, as they exist in undirected graphs, between effective resistances and control-theoretic properties, including robustness of linear consensus to noise \cite{Young10,Young11}, and node certainty in networks of stochastic decision-makers \cite{Poulakakis12}. Further applications of this concept to directed graphs should be possible in formation control \cite{Barooah06}, distributed estimation \cite{Barooah07, Barooah08} and optimal leader selection in networked control systems \cite{Patterson10, Clark11, Fardad11}.

Effective resistances have proved to be important in the study of networked systems because they relate global network properties to the individual connections between nodes, and they relate local network changes (e.g. the addition or deletion of an edge, or the change of an edge weight) to global properties without the need to re-compute these properties for the entire network (since only resistances that depend on the edge in question will change). Accordingly, the concept of effective resistance for directed graphs will be most useful if the resistance of any given connection within a graph can be computed, and if it is understood how to combine resistances from multiple connections. Computation and combination of resistances are possible for undirected graphs using the familiar rules for combining resistors in series and parallel. 

In this paper, we address the problems of computing and combining effective resistances for directed graphs. In Section \ref{sec:back} we review our definition of effective resistance for directed graphs from \cite{Young13I}. In Section \ref{sec:equalres} we develop some theory to identify directed graphs that have the same resistances as an equivalent undirected graph. We use these results in Section \ref{sec:direct} to recover the series-resistance formula for nodes connected by one directed path and the parallel-resistance formula for nodes connected by two directed paths in the form of a directed cycle. In Section \ref{sec:indirect} we examine nodes connected by a directed tree and derive a resistance formula that has no analogue from undirected graphs.

\section{Background and notation}\label{sec:back}
We present below some basic definitions of directed graph theory, as well as our definition of effective resistance. For more detail, the reader is referred to the companion paper \cite{Young13I}.  

A \emph{graph} $\mathcal{G}$ consists of the triple $\left(\mathcal{V}, \mathcal{E}, A \right)$, where $\mathcal{V} = \left\{1, 2, \ldots, N \right\}$ is the set of nodes, $\mathcal{E} \subseteq \mathcal{V}\times\mathcal{V}$ is the set of edges and $A \in \mathbb{R}^{N\times N}$ is a weighted adjacency matrix with non-negative entries $a_{i,j}$. Each $a_{i,j}$ will be positive if and only if $\left( i,j \right) \in \mathcal{E}$, otherwise $a_{i,j} = 0$. The graph $\mathcal{G}$ is said to be \emph{undirected} if $\left(i,j\right) \in \mathcal{E}$ implies $\left(j,i\right) \in \mathcal{E}$ and $a_{i,j} = a_{j,i}$. Thus, a graph will be undirected if and only if its adjacency matrix is symmetric.

The \emph{out-degree} of node $k$ is defined as $d_k^{\text{\emph{out}}} = \sum_{j=1}^N{a_{k,j}}$. $\mathcal{G}$ has an associated \emph{Laplacian} matrix $L$, defined by $L = D - A$, where $D$ is the diagonal matrix of node out-degrees.

A \emph{connection} in $\mathcal{G}$ between nodes $k$ and $j$ consists of two paths, one starting at $k$ and the other at $j$ and which both terminate at the same node. A \emph{direct connection} between nodes $k$ and $j$ is a connection in which one path is trivial (i.e. either only node $k$ or only node $j$) - thus a direct connection is equivalent to a path. Conversely, an \emph{indirect connection} is one in which the terminal node of the two paths is neither node $k$ nor node $j$.

The graph $\mathcal{G}$ is \emph{connected} if it contains a globally reachable node. Equivalently, $\mathcal{G}$ is connected if and only if a connection exists between any pair of nodes.

A \emph{connection subgraph} between nodes $k$ and $j$ in the graph $\mathcal{G}$ is a maximal connected subgraph of $\mathcal{G}$ in which every node and edge form part of a connection between nodes $k$ and $j$ in $\mathcal{G}$. If only one connection subgraph exists in $\mathcal{G}$ between nodes $k$ and $j$, it is referred to as \emph{the} connection subgraph and is denoted by $\mathcal{C}_\mathcal{G}(k,j)$. 

Let $Q \in \mathbb{R}^{(N-1)\times N}$ be a matrix that satisfies
\begin{equation}\label{eqn:propq}
Q\mathbf{1}_N = \mathbf{0}, \; QQ^T = I_{N-1} \text{ and } Q^TQ = \Pi.
\end{equation}

Using $Q$, we can compute the reduced Laplacian matrix for any graph as
\begin{equation}\label{eqn:lbar}
\overline{L} = QLQ^T,
\end{equation}
and then for connected graphs we can find the unique solution $\Sigma$ to the Lyapunov equation
\begin{equation}\label{eqn:lyap}
\overline{L}\Sigma + \Sigma \overline{L}^T = I_{N-1}.
\end{equation}
If we let
\begin{equation}\label{eqn:xdef}
X \mathrel{\mathop :}= 2Q^T \Sigma Q,
\end{equation}
the resistance between two nodes in a graph can be computed as
\begin{equation}\label{eqn:dirres}
r_{k,j} \!=\! \left(\mathbf{e}_N^{(k)} \!-\! \mathbf{e}_N^{(j)}\right)^{\!T}\!\!\! X \!\left(\mathbf{e}_N^{(k)} \!-\! \mathbf{e}_N^{(j)}\right) \!=\! x_{k,k} + x_{j,j} - 2x_{k,j}.
\end{equation}

Note that Definition \ref{P1:def:generalres} in the companion paper \cite{Young13I} extends effective resistance computations to disconnected graphs as well.

In some of the following results, we make use of \emph{binomial coefficients}, defined as
\begin{equation}\label{eqn:bincoef}
\binom{n}{k} = \frac{n!}{k! \left(n-k\right)!} \; n,k \in \mathbb{Z}, \; 0 \leq k \leq n.
\end{equation}

\section{Directed and undirected graphs with equal effective resistances}\label{sec:equalres}
In this section we prove Proposition \ref{prop:DAP}, which provides sufficient conditions for the resistances in a directed graph to be the same as the resistances in an equivalent undirected graph. The proof relies on two lemmas, which we prove first.

Recall that a \emph{permutation matrix} is a square matrix containing precisely one entry of $1$ in each row and column with every other entry being $0$.

\begin{lemma}\label{lem:Pprop}
Let $P$ be a permutation matrix. Then $P$ has the following properties
\begin{enumerate}[(i) ]
\item \begin{bulletequation}\label{eqn:Porthog}P^{-1}  = P^T,\end{bulletequation}
\item \begin{bulletequation}\label{eqn:PPi}P\Pi = \Pi P \text{, and}\end{bulletequation}
\item \begin{bulletequation}\label{eqn:PIPi}\left(P - I\right)\Pi = \Pi \left(P - I\right) = P - I. \end{bulletequation}
\end{enumerate}
\end{lemma}

\begin{IEEEproof}
\begin{enumerate}[(i) ]
\item This follows from the fact that the rows (or columns) of $P$ form an orthonormal set. See, e.g. Theorem 2.1.4 in \cite{Horn85}.

\item Since $P$ contains precisely one $1$ in each row and column, $P\mathbf{1}_N = \mathbf{1}_N$ and $\mathbf{1}_N^T = \mathbf{1}_N^TP$. Thus $P\Pi = P - \frac{1}{N}P\mathbf{1}_N\mathbf{1}_N^T = P - \frac{1}{N}\mathbf{1}_N\mathbf{1}_N^TP = \Pi P$.

\item The first part follows from \refe{eqn:PPi}, and again using the fact that $P\mathbf{1}_N = \mathbf{1}_N$ and $\mathbf{1}_N^T = \mathbf{1}_N^TP$, we have $(P-I)\Pi = P - I - \frac{1}{N}P\mathbf{1}_N\mathbf{1}_N^T + \frac{1}{N}\mathbf{1}_N\mathbf{1}_N^T = P - I$. \hfill \IEEEQEDhere

\end{enumerate}
\end{IEEEproof}

Since $P^T$ also satisfies the requirements of a permutation matrix, the results of Lemma \ref{lem:Pprop} apply to $P^T$ as well (this can also be seen by simply transposing equations \refe{eqn:Porthog}, \refe{eqn:PPi} and \refe{eqn:PIPi}).

The following lemma is required to prove Proposition \ref{prop:DAP}.
\begin{lemma}\label{lem:APdiag}
Let $A$ be a square matrix and $P$ be a permutation matrix of the same dimension as $A$. Suppose that $AP$ is diagonal. Then
\begin{enumerate}[(i) ]
\item $PA$ is also diagonal, \label{lempart:diagtoo}
\item $A\left(P - I\right) + A^T\left(P^T - I\right)$ is symmetric, that is
\begin{multline}\label{eqn:APIsymm}
A\left(P - I\right) + A^T\left(P^T - I\right) =\\ \left(P - I\right)A + \left(P^T - I\right)A^T \text{, and}
\end{multline}
\item \begin{bulletequation}\label{eqn:PAbars}\overline{(P - I)}^T \overline{A}^T \overline{A}\,\overline{(P - I)} = \overline{(P - I)} \,\overline{A} \, \overline{A}^T\overline{(P - I)}^T.\end{bulletequation}
\end{enumerate}
\end{lemma}

\begin{IEEEproof}
\begin{enumerate}[(i) ]
\item Let $D \mathrel{\mathop :}= AP$, which is diagonal by assumption. Then, by Lemma \ref{lem:Pprop}, we can see that $A = D P^T$. Thus $PA = P D P^T$, which implies that $PA$ is formed by permuting the rows and columns of a diagonal matrix, and is therefore diagonal.

\item Since $AP$ and $PA$ are diagonal, they are both symmetric. Thus $A^T P^T$ (${} = \left(PA\right)^T = PA$) is symmetric too. Since $\left(-A - A^T\right)$ is also symmetric, the result follows.

\item First we note that as $A P$ is diagonal, it is symmetric and commutes with its transpose (i.e. itself). Thus $P^T A^T A P = A P P^T A^T = A A^T$ (by \refe{eqn:Porthog}). Similarly, by part (\ref{lempart:diagtoo}), $P A$ is also diagonal and so it too is symmetric and commutes with its transpose. Hence $P A A^T P^T = A^T P^T P A = A^T A$ (by \refe{eqn:Porthog}). Using these facts, we can observe that $\left(P^T\! -\! I\right)\! A^T A\! \left(P\! -\! I\right) = \left(P\! -\! I\right)\!A A^T\! \left(P^T\! -\! I\right)$. Now, adding $\left(P\! -\! I\right)\!A^2\!\left(P\! -\! I\right)$ to both sides gives us $\left[\left(P - I\right)\!A + \left(P^T - I\right)\!A^T\right]\! A \left(P - I\right) = \left(P - I\right) A \left[A\!\left(P - I\right) + A^T\! \left(P^T - I\right)\right]$. But we can use \refe{eqn:APIsymm} to write this as
\begin{multline}
A\left(P - I\right)A\left(P - I\right) + A^T\left(P^T - I\right)A\left(P - I\right) = \\ 
\left(P \!-\! I\right) A\left(P \!-\! I\right) A + \left(P \!-\! I\right) A\left(P^T \!-\! I\right) A^T. \label{eqn:APproof1}
\end{multline}

Now, by \refe{eqn:PIPi}, we can pre- or post-multiply any factor of $\left(P - I\right)$ or $\left(P^T - I\right)$ by $\Pi$ without changing the matrix. Therefore, we can subtract $\left(P - I\right)A\Pi A\left(P - I\right)$ from both sides of \refe{eqn:APproof1}, obtain a common factor of $\Pi A \left(P - I\right)$ on the left hand side and $\left(P - I\right)A\Pi$ on the right hand side, then use \refe{eqn:APIsymm} to obtain 
\[
\left(P^T \!-\! I\right)A^T\Pi A \left(P \!-\! I\right) = \left(P \!-\! I\right) A \Pi A^T \left(P^T \!-\! I\right),
\]
which is equivalent to (using \refe{eqn:PIPi} again)
\begin{multline*}
\left(P^T - I\right)\Pi A^T\Pi A \Pi \left(P - I\right) = \\
\left(P - I\right)\Pi A \Pi A^T \Pi \left(P^T - I\right).
\end{multline*}
Finally, pre-multiplying by $Q$ and post-multiplying by $Q^T$ gives us our desired result.\hfill \IEEEQEDhere
\end{enumerate}
\end{IEEEproof}

The following proposition provides sufficient conditions for a directed graph to have the same resistances between any pair of nodes as an equivalent undirected graph. Although the assumption for the proposition may seem relatively general, it is straightforward to show that this can only apply to directed path and cycle graphs.

\begin{proposition}\label{prop:DAP}
Suppose $\mathcal{G} = \left(\mathcal{V},\mathcal{E},A\right)$ is a connected (directed) graph with matrix of node out-degrees $D$. Furthermore, suppose there is a permutation matrix $P$ such that $D = AP$. Let $\mathcal{G}_u = \left(\mathcal{V}_u,\mathcal{E}_u,A_u\right)$ be the undirected graph with $\mathcal{V}_u = \mathcal{V}$, $\mathcal{E}_u$ such that $(i,j) \in \mathcal{E} \Rightarrow (i,j) \text{ and } (j,i) \in \mathcal{E}_u$, and $A_u = \frac{1}{2}\left(A + A^T\right)$. Then the effective resistance between two nodes in $\mathcal{G}$ is equal to the effective resistance between the same two nodes in $\mathcal{G}_u$.
\end{proposition}

\begin{IEEEproof}
First we note that $\mathcal{G}_u$ must be an undirected graph since its adjacency matrix is symmetric. The Laplacian matrix $L$ of $\mathcal{G}$ is given by $L = D - A = A\left(P - I\right)$. Thus $\overline{L}$ is given by $\overline{L}	= QA\left(P - I\right)Q^T$, which can be rewritten (using \refe{eqn:PIPi}) as $\overline{L} = \overline{A} \, \overline{\left(P - I\right)}$. Furthermore, since $\mathcal{G}$ is connected, $\overline{L}$ is invertible \cite{Young10}.

Next, we claim that the Laplacian matrix of $\mathcal{G}_u$ is given by $L_u\! \mathrel{\mathop :}=\! \frac{1}{2}\!\left[A\left(P - I\right) + A^T\left(P^T - I\right)\right]$. To see this, we first note that we can rewrite $L_u$ as $L_u = \frac{1}{2}\left(D + A^TP^T\right) - A_u$. But by part (\ref{lempart:diagtoo}) of Lemma \ref{lem:APdiag}, $PA$ is diagonal and therefore so is $A^TP^T$. Hence $D_u \mathrel{\mathop :}= \dfrac{1}{2}\left(D + A^TP^T\right)$ is a diagonal matrix. Furthermore, $L_u \mathbf{1}_N = \frac{1}{2}\left[A\left(P\mathbf{1}_N - \mathbf{1}_N\right) + A^T\left(P^T\mathbf{1}_N - \mathbf{1}_N\right)\right]$. But since $P$ is a permutation matrix, $P\mathbf{1}_N = \mathbf{1}_N$ and $P^T\mathbf{1}_N = \mathbf{1}_N$, and so $L_u \mathbf{1}_N = \mathbf{0}$. Therefore, $L_u$ is equal to a diagonal matrix minus $A_u$ and has zero row sums. Hence $D_u$ must be the diagonal matrix of the row sums of $A_u$, i.e. the matrix of out-degrees of $\mathcal{G}_u$.

Since $\mathcal{G}_u$ contains every edge in $\mathcal{G}$ (in addition to the reversal of each edge) and $\mathcal{G}$ is connected, $\mathcal{G}_u$ must also be connected. Thus $\Sigma_{u} = \frac{1}{2}\overline{L}_u^{-1}$ is the solution to the Lyapunov equation \refe{eqn:lyap} for $\mathcal{G}_u$. Using our expression for $L_u$ and \refe{eqn:PIPi}, we can write $\Sigma_{u} = \left[\overline{A}\,\overline{\left(P - I\right)} + \overline{A}^T \, \overline{\left(P - I\right)}^T\right]^{-1}$. Since $\Sigma_u$ is symmetric, we can also write $\Sigma_{u} = \left[\overline{\left(P - I\right)}\,\overline{A} + \overline{\left(P - I\right)}^T\,\overline{A}^T\right]^{-1}$.

Now, when we substitute $\Sigma_{u}$ into the left hand side of equation \refe{eqn:lyap} for $\mathcal{G}$, we obtain
\begin{multline*}
\overline{L}\Sigma_{u} + \Sigma_{u}\overline{L}^T = \left[I + \overline{A}^T\,\overline{\left(P - I\right)}^T\left(\overline{A}\,\overline{\left(P - I\right)}\right)^{-1}\right]^{-1} \\
{} + \left[I + \left(\overline{\left(P - I\right)}^T\,\overline{A}^T\right)^{-1}\overline{\left(P - I\right)}\,\overline{A}\right]^{-1}.
\end{multline*}
Using the Matrix Inversion Lemma \cite{Woodbury50} applied to the first term, we can rewrite this as
\begin{multline}\label{eqn:SigmaPAbars}
\overline{L}\Sigma_{u} + \Sigma_{u}\overline{L}^T = I - \left[I + \overline{A}\,\overline{\left(P - I\right)}\left(\overline{A}^T\overline{\left(P - I\right)}^T\right)^{-1}\right]^{-1} \\
{} + \left[I + \left(\overline{\left(P - I\right)}^T\overline{A}^T\right)^{-1}\overline{\left(P - I\right)}\,\overline{A}\right]^{-1}\!\!\!.
\end{multline}

But by \refe{eqn:PAbars}, $\overline{(P - I)}^T \overline{A}^T \overline{A}\,\overline{(P - I)} = \overline{(P - I)} \,\overline{A} \, \overline{A}^T\overline{(P - I)}^T$, so $\overline{A}\,\overline{(P - I)}\left(\overline{A}^T\overline{(P - I)}^T\right)^{-1} = \left(\overline{(P - I)}^T \overline{A}^T\right)^{-1}\overline{(P - I)} \,\overline{A}$ and the final two terms in \refe{eqn:SigmaPAbars} are equal (with opposite signs). Thus $\overline{L}\Sigma_{u} + \Sigma_{u}\overline{L}^T = I$, and so $\Sigma_{u}$ solves \refe{eqn:lyap} for $\mathcal{G}$. This implies $\Sigma = \Sigma_{u}$, $X = X_u$ and $r_{k,j} = r_{u \; k,j}$ for all nodes $k$ and $j$.
\end{IEEEproof}

\section{Effective resistances from direct connections}\label{sec:direct}

In this section we compute the resistance in directed graphs between a pair of nodes that are only connected through a single direct connection, or two direct connections in opposite directions (i.e. the connection subgraph consists of either a directed path or a directed cycle). These two scenarios are analogous (in undirected graphs) to combining multiple resistances in series and combining two resistances in parallel. At present, we do not have general rules for combining resistances from multiple direct connections.

The most basic connection is a single directed edge. Intuitively, since an undirected edge with a given weight is equivalent to two directed edges (in opposite directions) with the same weight, one would expect that the resistance of a directed edge should be twice that of an undirected edge with the same weight. The following lemma shows that this is indeed true.

\begin{lemma}\label{lem:edgeres}
If $\mathcal{C}_\mathcal{G}(k,j)$ consists of a single directed edge from node $k$ to node $j$ with weight $a_{k,j}$, then
\begin{equation}\label{eqn:edgeres}
r_{k,j} = \frac{2}{a_{k,j}}.
\end{equation}
\end{lemma}

\begin{IEEEproof}
If we take node $j$ to be the first node in $\mathcal{C}_\mathcal{G}(k,j)$ and node $k$ to be the second, then $\mathcal{C}_\mathcal{G}(k,j)$ has Laplacian matrix $L = \begin{bmatrix}0 & 0\\-a_{k,j} & a_{k,j}\end{bmatrix}$. In this case, there is only one matrix $Q$ (up to a choice of sign) which satisfies \refe{eqn:propq}, namely $Q = \begin{bmatrix}\frac{1}{\sqrt{2}} & -\frac{1}{\sqrt{2}}\end{bmatrix}$. Then we have $\overline{L} = Q L Q^T = a_{k,j}$, and hence $\Sigma = \frac{1}{2a_{k,j}}$. Thus $X = 2Q^T\Sigma Q = \begin{bmatrix}\frac{1}{2a_{k,j}} & -\frac{1}{2a_{k,j}}\\-\frac{1}{2a_{k,j}} & \frac{1}{2a_{k,j}}\end{bmatrix}$, and finally,
\[
r_{k,j} = \left(\mathbf{e}_2^{(1)} - \mathbf{e}_2^{(2)}\right)^T X \left(\mathbf{e}_2^{(1)} - \mathbf{e}_2^{(2)}\right) = \frac{2}{a_{k,j}}. \IEEEQEDhereeqn
\]
\end{IEEEproof}

As a result of Lemma \ref{lem:edgeres}, when we refer to the effective resistance of a single (directed) edge, we mean twice the inverse of the edge weight. Our next two results extend to some directed graphs the familiar rules from undirected graphs for combining resistances in series and parallel. These cover the cases when a pair of nodes is connected only by either a directed path or cycle.

\begin{theorem}\label{theo:pathres}
Suppose $\mathcal{C}_\mathcal{G}(k,j)$ consists of a single directed path. Then $r_{k,j}$ is given by the sum of the resistances of each edge in the path between the two nodes (where the resistance of each edge is computed as in Lemma \ref{lem:edgeres}).
\end{theorem}

\begin{IEEEproof}
Suppose we label the nodes in $\mathcal{C}_\mathcal{G}(k,j)$ from $1$ to $N$ in the order in which they appear along the path, starting with the root and moving in the direction opposite the edges. Then we can write the adjacency matrix of $\mathcal{C}_\mathcal{G}(k,j)$ as $A = \text{diag}^{(-1)}\!\left(\begin{bmatrix}a_1 & a_2 & \cdots & a_{N-1}\end{bmatrix}\right)$, and the matrix of node out-degrees as $D = \text{diag}\!\left(\begin{bmatrix}0 & a_1 & \cdots & a_{N-1}\end{bmatrix}\right)$.

If we let $P$ be the permutation matrix containing ones above the main diagonal and in the lower left corner, we can observe that $D = AP$. Therefore, by Proposition \ref{prop:DAP}, the resistance between any two nodes in $\mathcal{C}_\mathcal{G}(k,j)$ is equal to the resistance between the same two nodes in an undirected graph with adjacency matrix $A_u = \frac{1}{2}\left(A + A^T\right)$.

Now, $A_u$ is the adjacency matrix of an undirected path, with weights of $\dfrac{1}{2}a_i$ on each edge. But the resistance of an edge in an undirected graph is the inverse of the edge weight and so each edge has resistance $\dfrac{2}{a_i}$. Thus the edge resistances in this undirected graph match those in the original directed path graph (computed according to Lemma \ref{lem:edgeres}). Furthermore, the resistance between two nodes connected by an undirected path is simply the sum of the resistances of the edges between them. Thus the same is true for two nodes connected by a directed path.
\end{IEEEproof}

\begin{theorem}\label{theo:cycleres}
Suppose $\mathcal{C}_\mathcal{G}(k,j)$ consists of a single directed cycle. Then $r_{k,j}$ is given by the inverse of the sum of the inverses of the resistances of each path connecting nodes $k$ and $j$ (where the resistance of each path is computed as in Theorem \ref{theo:pathres}).
\end{theorem}

\begin{IEEEproof}
Suppose we label the nodes in $\mathcal{C}_\mathcal{G}(k,j)$ from $1$ to $N$ in the reverse of the order in which they appear around the cycle, starting with any node. Then we can write the adjacency matrix of $\mathcal{C}_\mathcal{G}(k,j)$ as $A = \text{diag}^{(N-1)}\!\left(\begin{bmatrix}a_1\end{bmatrix}\right) + \text{diag}^{(-1)}\!\left(\begin{bmatrix}a_2 & a_3 & \cdots & a_N\end{bmatrix}\right)$ and the matrix of node out-degrees as $D = \text{diag}\!\left(\begin{bmatrix}a_1 & a_2 & \cdots & a_N\end{bmatrix}\right)$.

If we let $P$ be the permutation matrix containing ones above the main diagonal and in the lower left corner, we can observe that $D = AP$. Therefore, by Proposition \ref{prop:DAP}, the resistance between any two nodes in $\mathcal{C}_\mathcal{G}(k,j)$ is equal to the resistance between the same two nodes in an undirected graph with adjacency matrix $A_u = \frac{1}{2}\left(A + A^T\right)$.

Now, $A_u$ is the adjacency matrix of an undirected cycle, with weights of $\dfrac{1}{2}a_i$ on each edge. But the resistance of an edge in an undirected graph is the inverse of the edge weight, so each edge has resistance $\dfrac{2}{a_i}$. Thus the edge resistances in this undirected graph match those in the original directed cycle graph (computed according to Lemma \ref{lem:edgeres}). Furthermore, the resistance between nodes $k$ and $j$ connected by an undirected cycle is given by
\[
r_{u \, k,j} = \frac{1}{\frac{1}{r_1} + \frac{1}{r_2}},
\]
where $r_1$ is the resistance of one path between nodes $k$ and $j$ and $r_2$ is the resistance of the other path. Thus the same is true for two nodes connected by a directed cycle, where (by Theorem \ref{theo:pathres}) $r_1$ and $r_2$ are equal to the resistances of the two directed paths between nodes $k$ and $j$.
\end{IEEEproof}

\section{Effective resistances from indirect connections}\label{sec:indirect}
Lemma \ref{lem:edgeres} and Theorems \ref{theo:pathres} and \ref{theo:cycleres} suggest a very intuitive interpretation of effective resistance for directed graphs. A directed edge can be thought of as ``half'' of an undirected edge - either by noting that a directed edge allows half of the interaction to take place that occurs through an undirected edge, or by viewing an undirected edge as consisting of two directed edges with equal weights but in opposite directions. Thus, the resistance of a directed edge is twice the resistance of an undirected edge with the same weight. Then, in path and cycle graphs, resistances combine in exactly the ways (i.e. in series and in parallel) we are used to. However, connections in directed graphs can be more complicated than these. In particular, two nodes in a directed graph may be connected even if neither node is reachable from the other. This will occur when the only connections between the nodes consist of two non-zero length paths which meet at a distinct node. In Theorem \ref{theo:treeres} we prove an explicit expression for resistances in the case when $\mathcal{C}_\mathcal{G}(k,j)$ is a directed tree with unit edge weights. Before doing so we prove two lemmas on the correspondence between resistances and the matrix $X$ from \refe{eqn:xdef}, and two lemmas on the resistance between two leaves in a directed tree. We also rely on the finite series expressions given and proved in Appendix \ref{apndx:finser}.

\begin{lemma}\label{lem:Xres}
There is a one-to-one relationship between the effective resistances between nodes in a graph and the entries of the matrix $X$ from \refe{eqn:xdef}. In particular,
\begin{equation}\label{eqn:xtor}
r_{k,j} = x_{k,k} + x_{j,j} - 2x_{k,j} \text{, and}
\end{equation}
\begin{equation}\label{eqn:rtox}
x_{k,j} = \frac{1}{2N}\sum_{i = 1}^N{r_{k,i}} + \frac{1}{2N}\sum_{i = 1}^N{r_{j,i}} - \frac{1}{N^2}\sum_{i = 1}^{N-1}{\sum_{\ell = i+1}^N{r_{i,\ell}}} - \frac{1}{2}r_{k,j}.
\end{equation}
\end{lemma}

\begin{IEEEproof}
\refe{eqn:xtor} is simply the definition of $r_{k,j}$. To derive \refe{eqn:rtox}, we first note that from \refe{eqn:propq} and \refe{eqn:xdef}, $X$ has the property that $X\mathbf{1}_N = \mathbf{0}$ and $\mathbf{1}_N^T X	= \mathbf{0}^T$. That is, $X$ has zero row- and column-sums.

Now, using \refe{eqn:xtor}, we can write $r_{k,i} = x_{k,k} + x_{i,i} - 2x_{k,i}$ for any $1 \leq i \leq N$. Then, by summing this equation over $i$, we obtain $\sum_{i = 1}^N{r_{k,i}} = Nx_{k,k} + \text{tr}\left(X\right)$ (since $X$ has zero row-sums). Next, by summing again over $k$, we find that $\text{tr}\left(X\right) = \frac{1}{2N} \sum_{k = 1}^N{\sum_{i = 1}^N{r_{k,i}}}$. But $r_{i,i} = 0 \; \forall i$ and $r_{i,k} = r_{k,i}$ (by Theorem \ref{P1:theo:metric} in the companion paper \cite{Young13I}). Thus we can say that
\begin{equation}\label{eqn:traceX}
\text{tr}\left(X\right) = \frac{1}{N} \sum_{i = 1}^{N-1}{\sum_{\ell = i+1}^N{r_{i,\ell}}}.
\end{equation}

Combining \refe{eqn:traceX} with our expression for $\sum_{i = 1}^N{r_{k,i}}$ gives us
\begin{equation}\label{eqn:xkk}
x_{k,k} = \frac{1}{N}\sum_{i = 1}^N{r_{k,i}} - \frac{1}{N^2}\sum_{i = 1}^{N-1}{\sum_{\ell = i+1}^N{r_{i,\ell}}}.
\end{equation}
Substituting the expression from \refe{eqn:xkk} for $x_{k,k}$ and $x_{j,j}$ in \refe{eqn:xtor} produces \refe{eqn:rtox}.
\end{IEEEproof}

\begin{lemma}\label{lem:Xpath}
Suppose $\mathcal{G}$ is a directed path with unit edge weights containing $N$ nodes, in which the nodes are labelled from $1$ to $N$ in the order in which they appear along the path, starting with the root. Let $X$ be the corresponding matrix from \refe{eqn:xdef}. Then the entries of $X$ are given by
\begin{multline}\label{eqn:xpath}
x_{k,j} \!=\! \frac{2N^2 \!+ 3N \!+ 1 + 3k^2 \!+ 3j^2 \!- 3\left(N\!+\!1\right)k - 3\left(N\!+\!1\right)j}{3N} \\
{} - \left|k - j\right|.
\end{multline}
\end{lemma}

\begin{IEEEproof}
Suppose $k, j \, \in \, \left\{1,2,\ldots,N\right\}$. Then by Theorem \ref{theo:pathres}, we know that the resistance between nodes $k$ and $j$ in our directed path is equal to $2$ (the resistance of each edge) times the number of edges between them. Since the nodes are labelled in order along the path, this gives us $r_{k,j} = 2\left|k - j\right|$. Therefore, from Lemma \ref{lem:Xres}, we know that
\begin{multline}\label{eqn:xpathinit}
x_{k,j}	= \frac{1}{N}\sum_{i = 1}^N{\left|k - i\right|} + \frac{1}{N}\sum_{i = 1}^N{\left|j - i\right|} - \frac{2}{N^2}\sum_{i = 1}^{N-1}{\sum_{\ell = i+1}^N{\!\!\left|i - \ell\right|}} \\
{} - \left|k - j\right|.
\end{multline}

We now proceed by examining each summation in turn. The first sum can be broken into two parts and then simplified using \refe{eqn:sumint} to obtain $\sum_{i=1}^N{\left|k - i\right|}	= \frac{2k^2 - 2\left(N+1\right)k + \left(N+1\right)N}{2}$. By replacing $k$ with $j$ in the previous expression, we observe that $\sum_{i=1}^N{\left|j - i\right|} = \frac{2j^2 - 2\left(N+1\right)j + \left(N+1\right)N}{2}$.

In the third sum in \refe{eqn:xpathinit}, we observe that $\ell > i$ for every term. Thus $\left|i - \ell\right| = \ell - i$, and we can use \refe{eqn:sumint} and \refe{eqn:sumintsq} to obtain $\sum_{i = 1}^{N-1}{\sum_{\ell = i+1}^N{\left|i - \ell\right|}} = \frac{\left(N^2 - 1\right)N}{6}$. Finally, \refe{eqn:xpath} follows from substituting our expressions for each summation into \refe{eqn:xpathinit}.
\end{IEEEproof}

The following results are needed to prove Theorem \ref{theo:treeres}. In them, we examine the resistance between the leaves of a tree containing two branches that meet at its root and with unit weights on every edge, $\mathcal{G}^\text{tree}_{n,m}$, as shown in Fig. \ref{fig:C_G_trees}\subref{fig:C_G_n_m}. The effective resistance between the two leaves of $\mathcal{G}^\text{tree}_{n,m}$ will be denoted by $r(n,m)$.

\begin{lemma}\label{lem:treeresn1}
The effective resistance between the two leaves of $\mathcal{G}^\text{tree}_{n,1}$ is given by
\begin{equation}\label{eqn:treeresn1}
r(n,1) = 2(n-1) + 2^{2-n}.
\end{equation}
\end{lemma}

\begin{IEEEproof}
The number of nodes in $\mathcal{G}^\text{tree}_{n,1}$ is $N = n+2$. Let us label the nodes in $\mathcal{G}^\text{tree}_{n,1}$ from $1$ to $n+1$, in the reverse order of the edges, along the branch of length $n$, starting with the root (thus the leaf of this branch is node $n+1$). Then the other leaf (with an edge connecting it to the root) will be node $N = n+2$. Thus the resistance we seek to find is $r(n,1) = r_{n+1,n+2}$.

Let $A_{N_p}^\text{path}$, $D_{N_p}^\text{path}$ and $L_{N_p}^\text{path}$ denote the adjacency matrix, matrix of out-degrees and Laplacian matrix of a directed path containing $N_p$ nodes and unit weights on every edge. Let the nodes in this path be labelled from $1$ to $N_p$ in the reverse of the order in which they appear, starting with the root. Thus $A_{N_p}^\text{path} = \text{diag}^{(-1)}\!\left(\mathbf{1}_{N_p-1}\right)$, $D_{N_p}^\text{path} = \text{diag}\!\left(\begin{bmatrix}0 & \mathbf{1}_{N_p-1}^T\end{bmatrix}\right)$ and $L_{N_p}^\text{path} = \text{diag}\!\left(\begin{bmatrix}0 & \mathbf{1}_{N_p-1}^T\end{bmatrix}\right) - \text{diag}^{(-1)}\!\left(\mathbf{1}_{N_p-1}\right)$. From these, we can observe that
\begin{equation}\label{eqn:oneL}
\mathbf{1}_{N_p}^T L_{N_p}^\text{path} = \mathbf{e}_{N_p}^{(N_p) T} - \mathbf{e}_{N_p}^{(1) T} \text{, and}
\end{equation}
\begin{equation}\label{eqn:eiL}
\mathbf{e}_{N_p}^{(i)T} L_{N_p}^\text{path} = \begin{cases}\mathbf{e}_{N_p}^{(i) T} - \mathbf{e}_{N_p}^{(i-1) T} &\text{if } 1 < i \leq N_p,\\\mathbf{0}^T &\text{if } i = 1.\end{cases}
\end{equation}
Next, we will let $Q_{N_p}$ be a $\left(N_p - 1\right)\times N_p$ matrix which satisfies \refe{eqn:propq}, and $\overline{L}_{N_p}^\text{path}$ and $\Sigma_{N_p}^\text{path}$ be derived from \refe{eqn:lbar} and \refe{eqn:lyap} using $L_{N_p}^\text{path}$ and $Q_{N_p}$. Let $X_{N_p}^\text{path} = 2Q_{N_p}^T \Sigma_{N_p}^\text{path} Q_{N_p}$, according to \refe{eqn:xdef}. Then, by Lemma \ref{lem:Xpath}, the entries of $X_{N_p}^\text{path}$ are given by \refe{eqn:xpath}.

Now, we can write the adjacency matrix, matrix of out-degrees and Laplacian matrix of $\mathcal{G}^\text{tree}_{n,1}$ as $A = \begin{bmatrix}A_{n+1}^\text{path} & \mathbf{0}\\\mathbf{e}_{n+1}^{(1) T} & 0\end{bmatrix}$, $D = \begin{bmatrix}D_{n+1}^\text{path} & \mathbf{0}\\\mathbf{0}^T & 1\end{bmatrix}$, and $L = \begin{bmatrix}L_{n+1}^\text{path} & \mathbf{0} \\ -\mathbf{e}_{n+1}^{(1) T} & 1\end{bmatrix}$. Next, let $Q = \begin{bmatrix}Q_{n+1} & \mathbf{0}\\\alpha\mathbf{1}_{n+1}^T & -\beta\end{bmatrix}$, where $\alpha = \frac{1}{\sqrt{(n+1)(n+2)}}$ and $\beta = \sqrt{\frac{n+1}{n+2}}$. Then $Q$ satisfies \refe{eqn:propq}. We can use \refe{eqn:lbar}, \refe{eqn:oneL} and the facts that $L_{n+1}^\text{path} \mathbf{1}_{n+1} = \mathbf{0}_{n+1}$ and $\beta\left(\alpha + \beta\right) = 1$ to express $\overline{L}$ as
\[
\overline{L} = \begin{bmatrix}\overline{L}_{n+1}^\text{path} & \mathbf{0}\\  \left(\beta - \alpha\right)\mathbf{e}_{n+1}^{(1) T}Q_{n+1}^T + \alpha\mathbf{e}_{n+1}^{(n+1) T}Q_{n+1}^T & 1\end{bmatrix}.
\]

In order to compute resistances in $\mathcal{G}^\text{tree}_{n,1}$, we must find the matrix $\Sigma$ which solves \refe{eqn:lyap}. Since we have partitioned $\overline{L}$ into a $2\times 2$ block matrix, we will do the same for $\Sigma$. Let $\Sigma = \begin{bmatrix} S & \mathbf{t}\\ \mathbf{t}^T & u\end{bmatrix}$, where $S \in \mathbb{R}^{n\times n}$ is a symmetric matrix, $\mathbf{t} \in \mathbb{R}^{n}$ and $u \in \mathbb{R}$. Then multiplying out the matrices in \refe{eqn:lyap} and equating blocks in this matrix equation gives us
\begin{align}
\overline{L}_{n+1}^\text{path} S + S\overline{L}_{n+1}^{\text{path}\, T}	&= I_n, \label{eqn:Tn1lyapS} \\
\overline{L}_{n+1}^\text{path} \mathbf{t} + \mathbf{t} + \left(\beta - \alpha\right) S Q_{n+1} \mathbf{e}_{n+1}^{(1)} \hspace{1.3cm} &\nonumber\\
{} + \alpha S Q_{n+1}\mathbf{e}_{n+1}^{(n+1)}						&= \mathbf{0} \text{, and} \label{eqn:Tn1lyapt} \\
2u + 2\left(\beta - \alpha\right)\mathbf{e}_{n+1}^{(1) T}Q_{n+1}^T\mathbf{t} \hspace{2.3cm} & \nonumber\\
{} + 2\alpha \mathbf{e}_{n+1}^{(n+1) T}Q_{n+1}^T\mathbf{t}			&= 1. \label{eqn:Tn1lyapu}
\end{align}

From \refe{eqn:Tn1lyapS}, it is clear that $S = \Sigma_{n+1}^\text{path}$. In addition, we can rewrite \refe{eqn:Tn1lyapu} as
\begin{equation}\label{eqn:Tn1uform}
u = \frac{1}{2} - \left(\beta - \alpha\right)\mathbf{e}_{n+1}^{(1) T}Q_{n+1}^T\mathbf{t} - \alpha \mathbf{e}_{n+1}^{(n+1) T}Q_{n+1}^T\mathbf{t}.
\end{equation}
In order to find a complete solution for $\Sigma$, we must solve \refe{eqn:Tn1lyapt} for $\mathbf{t}$. However, resistances are computed from $X$, which, if we let $\mathbf{v} \mathrel{\mathop :}= Q_{n+1}^T \mathbf{t} = \left[v_i\right]$ and use \refe{eqn:xdef}, can be written as
\begin{multline*}
X = \left[\begin{matrix} X^\text{path}_{n+1} + 2\alpha \mathbf{v}\mathbf{1}_{n+1}^T + 2\alpha\mathbf{1}_{n+1} \mathbf{v}^T + 2\alpha^2 u \mathbf{1}_{n+1}\mathbf{1}_{n+1}^T \\ -2\beta \mathbf{v}^T - 2\alpha\beta u \mathbf{1}_{n+1}^T \end{matrix}\right.\\
\left.\begin{matrix}  -2\beta \mathbf{v} - 2\alpha\beta u \mathbf{1}_{n+1}\\ 2\beta^2 u \end{matrix}\right].
\end{multline*}
Hence, to compute resistances in $\mathcal{G}^\text{tree}_{n,1}$, we need only compute $\mathbf{v}$, not $\mathbf{t}$. We can also note that as $X$ does not depend on our choice of $Q$ (by Lemma \ref{P1:lem:indofq} in the companion paper \cite{Young13I}), neither does $\mathbf{v}$. In fact, we can write \refe{eqn:Tn1uform} as $u = \frac{1}{2} + \left(\alpha - \beta\right)v_{1} - \alpha v_{n+1}$, and the resistance we seek as
\begin{multline}
r(n,1) = x^\text{path}_{n+1 \, n+1,n+1} + \left(\alpha +\beta\right)^2 + 2\left(\alpha +\beta\right)^2\left(\alpha - \beta\right)v_1 \\
{} + 2\left(\alpha + \beta\right)\left[2 - \alpha(\alpha + \beta)\right] v_{n+1}. \label{eqn:rn1fromv}
\end{multline}
Thus we only need to find $v_1$ and $v_{n+1}$ in order to compute $r(n,1)$.

Now, $v_i = \mathbf{e}_{n+1}^{(i) T}\mathbf{v} = \mathbf{e}_{n+1}^{(i) T} Q_{n+1}^T \mathbf{t}$. We will therefore proceed by left-multiplying \refe{eqn:Tn1lyapt} by $\mathbf{e}_{n+1}^{(i) T} Q_{n+1}^T$. Using the fact that $S = \Sigma_{n+1}^\text{path}$, we obtain
\begin{multline}\label{eqn:Tn1lyapv}
\mathbf{e}_{n+1}^{(i) T} Q_{n+1}^T Q_{n+1} L_{n+1}^\text{path} \mathbf{v} + v_i + \frac{\beta - \alpha}{2} \mathbf{e}_{n+1}^{(i) T} X_{n+1}^\text{path} \mathbf{e}_{n+1}^{(1)} \\
{} + \frac{\alpha}{2} \mathbf{e}_{n+1}^{(i) T} X_{n+1}^\text{path} \mathbf{e}_{n+1}^{(n+1)} = 0.
\end{multline}
But $\mathbf{e}_{n+1}^{(i) T} Q_{n+1}^T Q_{n+1} \!=\! \mathbf{e}_{n+1}^{(i) T} \left(I_{n+1} \!-\! \frac{1}{n+1}\mathbf{1}_{n+1}\mathbf{1}_{n+1}^T\right) \!=\! \mathbf{e}_{n+1}^{(i) T} \!-\! \frac{1}{n+1}\mathbf{1}_{n+1}^T$ by \refe{eqn:propq}, so by \refe{eqn:oneL} and \refe{eqn:eiL},
\[
\mathbf{e}_{n+1}^{(i) T} Q_{n+1}^T Q_{n+1}L_{n+1}^\text{path}\mathbf{v} = \begin{cases}\frac{1}{n+1}v_1 + v_i - v_{i-1} - \frac{1}{n+1}v_{n+1} \hspace{-1.4cm}& \\ &\hspace{-1.3cm} \text{if } 1 < i \leq n+1,\\ \frac{1}{n+1}v_1 - \frac{1}{n+1}v_{n+1} & \text{if } i = 1.\end{cases}
\]
Furthermore, using \refe{eqn:xpath}, we observe that
\begin{align*}
\mathbf{e}_{n+1}^{(i) T} X_{n+1}^\text{path} \mathbf{e}_{n+1}^{(1)}	&= x^\text{path}_{n+1 \, i,1} \\
													&= \frac{(2n+3)(n+2)}{3(n+1)} + \frac{i(i - 2n - 3)}{n+1} \text{, and}
\end{align*}
\[
\mathbf{e}_{n+1}^{(i) T} X_{n+1}^\text{path} \mathbf{e}_{n+1}^{(n+1)} = x^\text{path}_{n+1 \, i,n+1} = -\frac{n(n+2)}{3(n+1)} + \frac{i(i-1)}{n+1}.
\]

Substituting these expressions into \refe{eqn:Tn1lyapv} gives us
\begin{align}
v_i	&= \frac{1}{2}v_{i-1} - \frac{1}{2(n+1)}v_1 + \frac{1}{2(n+1)}v_{n+1} + f + g(i) \nonumber\\
	&\hspace{5cm} \text{ if } 1< i \leq n+1 \label{eqn:recurrvi} \\
v_1	&= \frac{1}{n+2}v_{n+1} + h,\label{eqn:recurrv1}
\end{align}
where $f = \frac{\left[(3\alpha -2\beta)n + 3(\alpha-\beta)\right](n+2)}{12(n+1)}$, $g(i) = \frac{i\left[-\beta i + 2(\beta-\alpha) n - 2\alpha + 3\beta\right]}{4(n+1)}$, and $h = \frac{\alpha n}{6} + \frac{(\alpha - \beta)n(2n+1)}{6(n+2)}$.

We can now recursively apply \refe{eqn:recurrvi} $n$ times, starting with $i = n+1$, to find
\begin{multline}\label{eqn:recurrvisoln}
v_{n+1} = 2^{-n}v_{1} - \frac{v_1}{n+1}\sum_{k = 1}^{n}{2^{-k}} + \frac{v_{n+1}}{n+1}\sum_{k = 1}^{n}{2^{-k}} + 2f\sum_{k = 1}^{n}{2^{-k}} \\
{} + 2\sum_{k = 1}^{n}{g(n+2-k)2^{-k}}.
\end{multline}
But we can write $g(n+2-k) = g_1 + g_2 k + g_3 k^2$, where $g_1 = \frac{(\beta - 2\alpha)(n+2)}{4}$, $g_2 = \frac{2\alpha n + 2\alpha + \beta}{4(n+1)}$, and $g_3 = \frac{-\beta}{4(n+1)}$. Therefore, by \refe{eqn:sumtwos}, \refe{eqn:suminttwos} and \refe{eqn:sumintsqtwos}, the sum involving $g(n+2-k)$ is
\begin{multline*}
\sum_{k = 1}^{n}{g(n+2-k)2^{-k}} = g_1\left(1 \!-\! 2^{-n}\right) + g_2\left[2 \!-\! \left(n+2\right)2^{-n}\right] \\
{} + g_3\left[6 - \left(n^2 + 4n + 6\right)2^{-n}\right].
\end{multline*}
Using this result and \refe{eqn:sumtwos}, \refe{eqn:recurrvisoln} becomes
\begin{multline}
\frac{n + 2^{-n}}{n+1} v_{n+1} = \frac{\left(n+2\right)2^{-n} - 1}{n+1}v_1 + 2f\left(1 - 2^{-n}\right) \\
{} + 2g_1\left(1 - 2^{-n}\right) + 2g_2\left[2 - \left(n+2\right)2^{-n}\right] \\
	{}  + 2g_3\left[6 - \left(n^2 + 4n + 6\right)2^{-n}\right].\label{eqn:recurrvn1}
\end{multline}

But now \refe{eqn:recurrvn1} and \refe{eqn:recurrv1} form a pair of linear equations in $v_1$ and $v_{n+1}$. Using the expressions for $f$, $g_1$, $g_2$, $g_3$ and $h$, along with the definitions of $\alpha$ and $\beta$, their solution is given by
\begin{equation}\label{eqn:vsols}
\begin{aligned}v_1	&= \frac{\alpha\left[-2n^2 + 5n - 6 + 6.2^{-n}\right]}{6} \text{ and}\\
v_{n+1} &= \frac{\alpha\left[n^2 + 2n - 12 + (6n+12) 2^{-n}\right]}{6}.\end{aligned}
\end{equation}
Finally, using \refe{eqn:xpath} and \refe{eqn:vsols} in \refe{eqn:rn1fromv}, along with the expressions for $\alpha$ and $\beta$, gives us \refe{eqn:treeresn1}.
\end{IEEEproof}

\begin{lemma}\label{lem:treeresnl1}
For positive integers $n$ and $\ell$, the resistance between the two leaves of $\mathcal{G}^\text{tree}_{n,\ell+1}$ satisfies the recurrence relation
\begin{multline}
r(n,\ell+1)	= \frac{-3n^2 + 3\ell^2 - 2n\ell -n + 5\ell + 2}{2(n+\ell+1)^2} \\
{} + \frac{\ell^2 + 2n\ell + 2n + 3\ell}{n+\ell+1}2^{-n} + \frac{n^2 + n + 2}{2(n+\ell+1)}2^{-\ell}\\
	{} + \frac{1}{4(n+\ell+1)} \sum_{k=1}^\ell {\textstyle\left(4 - \frac{2}{n+\ell+1} - 2^{k-\ell}\right)r(n,k)} \\
	{} - \frac{n+\ell+2}{2(n+\ell+1)} \sum_{k=1}^n {\textstyle\left(\frac{1}{n+\ell+1} - 2^{k-n}\right)r(k,\ell)}\\
	{} - \frac{1}{4(n+\ell+1)} \sum_{k=1}^n {\sum_{j=1}^\ell {\left(2^{1+k-n} - 2^{j-\ell}\right)r(k,j)}}. \label{eqn:treeresnl1}
\end{multline}
\end{lemma}

The proof of Lemma \ref{lem:treeresnl1} relies on similar ideas to the proof of Lemma \ref{lem:treeresn1}, and is given in Appendix \ref{apndx:treeres}. We now proceed to solve the recurrence relation given by Lemmas \ref{lem:treeresn1} and \ref{lem:treeresnl1} using several finite series results given in Appendix \ref{apndx:finser}.

\begin{theorem}\label{theo:treeres}
Suppose $\mathcal{C}_\mathcal{G}(k,j)$ consists of a directed tree with unit weights on every edge. Then $r_{k,j}$ is given by
\begin{equation}\label{eqn:treeres}
r_{k,j} = 2\left(n - m\right) + 2^{3 - n - m}\sum_{i = 1}^{\left\lfloor\frac{m+1}{2}\right\rfloor}{i \binom{n+m+2}{n+2i+1}},
\end{equation}
where $n$ is the length of the shortest path from node $k$ to a mutually reachable node and $m$ is the length of the shortest path from node $j$ to a mutually reachable node.
\end{theorem}

\begin{IEEEproof}
Since every node in $\mathcal{C}_\mathcal{G}(k,j)$ is reachable from either node $k$ or node $j$, if $\mathcal{C}_\mathcal{G}(k,j)$ is a tree then only nodes $k$ and $j$ can be leaves. But every tree has at least one leaf, so suppose that node $k$ is a leaf. If node $j$ is not a leaf, then $\mathcal{C}_\mathcal{G}(k,j)$ must be a directed path and node $j$ is the closest mutually reachable node to both nodes $k$ and $j$. Then $m = 0$, $n$ is the path length from $k$ to $j$ and \refe{eqn:treeres} reduces to $r_{k,j} = 2n$, which follows from Theorem \ref{theo:pathres}. Conversely, if node $j$ is a leaf but node $k$ is not, $\mathcal{C}_\mathcal{G}(k,j)$ must be a directed path and node $k$ is the closest mutually reachable node to both nodes $k$ and $j$. Then $n = 0$, $m$ is the path length from $j$ to $k$ and \refe{eqn:treeres} reduces to $r_{k,j} = -2m + 2^{3-m}\sum_{i = 1}^{\left\lfloor\frac{m+1}{2}\right\rfloor}{i \binom{m+2}{2i+1}}$. But by \refe{eqn:binoddsumi} and \refe{eqn:binoddsum} from Lemma \ref{lem:standardsums}, $\sum_{i = 1}^{\left\lfloor\frac{m+1}{2}\right\rfloor}{i \binom{m+2}{2i+1}} = m2^{m-1}$, and so \refe{eqn:treeres} becomes $r_{k,j} = 2m$, which follows from Theorem \ref{theo:pathres}.

Now, if both node $k$ and node $j$ are leaves, then $\mathcal{C}_\mathcal{G}(k,j)$ must be a directed tree with exactly two branches. Thus $\mathcal{C}_\mathcal{G}(k,j)$ must correspond to the tree shown in Fig. \ref{fig:C_G_trees}\subref{fig:C_G_tree} and $n$ and $m$ are the path lengths from nodes $k$ and $j$, respectively, to the point where the two branches meet. Furthermore, both $n$ and $m$ are at least $1$.

By Corollary \ref{P1:cor:removepath} from the companion paper \cite{Young13I}, we observe that the resistance between nodes $k$ and $j$ remains the same as we remove all the nodes of $\mathcal{C}_\mathcal{G}(k,j)$ from the root to the node where the two branches meet. Thus, $r_{k,j}$ can be computed as the resistance between the two leaves of the tree shown in Fig. \ref{fig:C_G_trees}\subref{fig:C_G_n_m}. Let this tree be called $\mathcal{G}^\text{tree}_{n,m}$, and since the only two parameters that define $\mathcal{G}^\text{tree}_{n,m}$ are $n$ and $m$, we can write $r_{k,j}$ as a function of $n$ and $m$ only. That is,
\[
r_{k,j} = \mathrel{\mathop :} r(n,m).
\]

\begin{figure}
\centering
\subfloat{{\footnotesize(a)}\includegraphics[width=2.5cm]{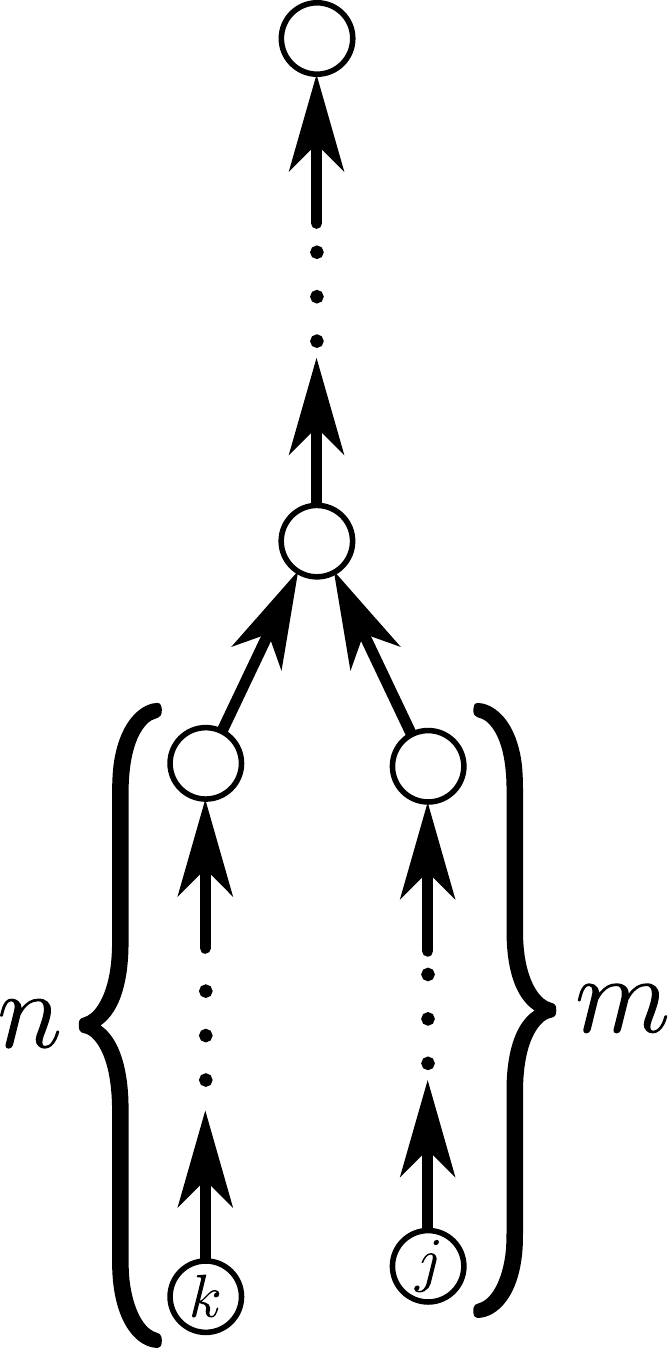}\label{fig:C_G_tree}}
\hspace{1cm}
\subfloat{{\footnotesize(b)}\includegraphics[width=2.5cm]{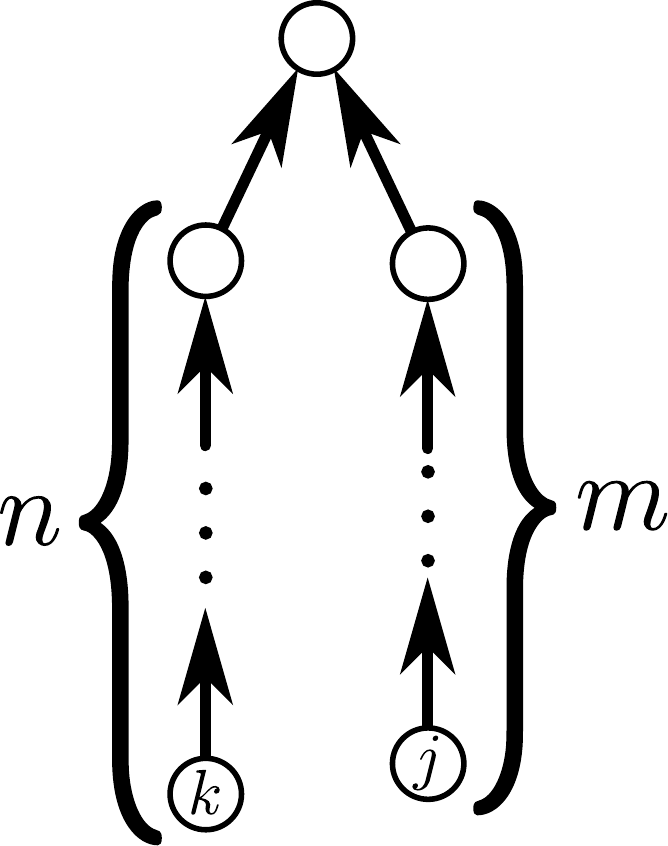}\label{fig:C_G_n_m}}
\caption{(a) The generic form of $\mathcal{C}_\mathcal{G}(k,j)$ when it is a directed tree with more than one leaf and unit weights on every edge. (b) A tree, $\mathcal{G}^\text{tree}_{n,m}$, in which $r_{k,j}$ is equal to its value in $\mathcal{C}_\mathcal{G}(k,j)$ when $\mathcal{C}_\mathcal{G}(k,j)$ is a directed tree as shown in part (a).}\vspace{-0.5cm}
\label{fig:C_G_trees}
\end{figure}

In order to compute $r(n,m)$, we will begin by considering the case where $m = 1$. Substituting $m = 1$ into \refe{eqn:treeres} gives $r(n,1) = 2(n-1) + 2^{2-n}$, which follows from Lemma \ref{lem:treeresn1}.

Now, suppose that \refe{eqn:treeres} holds for all $n > 0$ and all $m \in \left\{1, 2, \ldots, \ell\right\}$ (for some $\ell > 0$). Then $r_{k,j}$ for $m = \ell + 1$ can be computed using Lemma \ref{lem:treeresnl1}. In particular, all resistances in the right-hand side of \refe{eqn:treeresnl1} are given by \refe{eqn:treeres}. Therefore, we find that $r(n,\ell+1)$ matches the expression $s(n,\ell)$ given in Lemma \ref{lem:rsumsimp}. Therefore, $r(n,\ell+1)$ can be expressed in the form given in \refe{eqn:rsumsimp}.

Next, suppose that $\ell$ is odd. That is, $\ell = 2p + 1$ for some integer $p \geq 0$. Then \refe{eqn:rsumsimp} gives us
\begin{multline}
r(n, 2p+2) = 2\left(n - 2p - 2\right) + 2^{1-n-2p}\sum_{i=1}^{p+1}{i\binom{n+2p+4}{n+2i+1}} \\
{} + \frac{g(n,p)}{n+\ell+1},
\end{multline}
where $g(n,p)$ is given by \refe{eqn:gdef} in Lemma \ref{lem:oddexpression}. But by Lemma \ref{lem:oddexpression}, $g(n,p) = 0$ for any integers $n \geq 0$ and $p \geq 0$. Thus \refe{eqn:treeres} holds for $m = \ell + 1$.

Finally, suppose that $\ell$ is even. That is, $\ell = 2p$ for some integer $p > 0$. Then \refe{eqn:rsumsimp} gives us
\begin{multline*}
r(n, 2p+1)	= 2\left(n - 2p - 1\right) + 2^{2-n-2p}\sum_{i=1}^{p}{i\binom{n+2p+3}{n+2i+1}}  \\
	{} + \frac{4p^2 + 2np + 2n + 6p +2}{n+2p+1}2^{1-n-2p} + \frac{h(n,p)}{n+\ell+1},
\end{multline*}
where $h(n,p)$ is given by \refe{eqn:hdef} in Lemma \ref{lem:evenexpression}. But by Lemma \ref{lem:evenexpression}, $h(n,p) = 0$ for any integers $n \geq 0$ and $p \geq 0$. Thus,
\[
r(n, 2p+1)	= 2\left(n - 2p - 1\right) + 2^{2-n-2p}\sum_{i=1}^{p+1}{i\binom{n+2p+3}{n+2i+1}},
\]
and so \refe{eqn:treeres} holds for $m = \ell + 1$.

Therefore, by induction we have that \refe{eqn:treeres} also holds for all $n > 0$ and $m > 0$.
\end{IEEEproof}

Equation \refe{eqn:treeres} is a highly non-intuitive result, not least because on initial inspection it does not appear to be symmetric in $n$ and $m$ (although we know that it must be, by Theorem \ref{P1:theo:metric} in the companion paper). Therefore, it becomes easier to interpret \refe{eqn:treeres} if we reformulate it in terms of the shorter path length and the difference between the path lengths. Thus, if we let $n$ be the length of the longer path, that is, $n = m + d$ for some $d \geq 0$, \refe{eqn:treeres} becomes
\[
r_{k,j} = 2d + 2^{3-2m-d}\sum_{i=1}^{\left\lfloor\frac{m+1}{2}\right\rfloor}{i{\textstyle \binom{2m+d+2}{m+d+2i+1}}} =\mathrel{\mathop :} 2d + e(m,d).
\]
Then, using \refe{eqn:bincoef}, we can write
\begin{align*}
e(m,d+1)	&= 2^{3-2m-d}\!\!\sum_{i=1}^{\left\lfloor\frac{m+1}{2}\right\rfloor}{i\frac{2m+d+3}{2m\!+\!2d\!+\!4i\!+\!4}{\textstyle\binom{2m+d+2}{m+d+2i+1}}} \\
		&< \frac{2m+d+3}{2m+2d+4}e(m,d),
\end{align*}
and hence conclude that $\displaystyle\lim_{d\rightarrow\infty}e(m,d) = 0$. Thus, when the connection subgraph between two nodes is a directed tree, the resistance between them is twice the difference between the lengths of the paths connecting each node to their closest mutually reachable node, plus some ``excess'' that disappears as this difference becomes large. Conversely, the excess is significant when the path length difference is small, leading to a resistance that is greater than twice the difference.

\begin{figure}
\centering
\includegraphics[width=3cm]{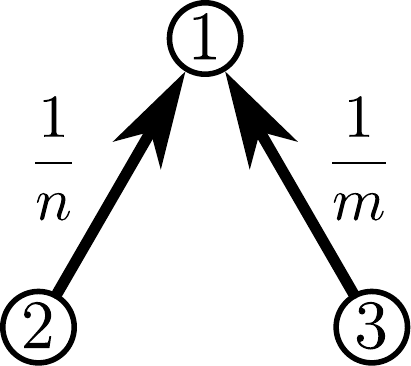}
\caption{A simple 3-node directed graph, $\mathcal{G}^\text{star}_3$, with resistances of $2 n$ and $2 m$ on each edge.}\vspace{-0.5cm}
\label{fig:indirectexample}
\end{figure}

One common approach to the analysis of resistive circuits is to replace a section of the network that connects to the rest through a single pair of nodes by a single resistor with an equivalent resistance. The simplest example of this is the replacement of a path with a single edge with equivalent resistance. If this principle were to extend to the calculation of effective resistance in directed graphs, then $r_{2,3}$ in $\mathcal{G}^\text{star}_3$ (as shown in Fig. \ref{fig:indirectexample}) would match the formula from Theorem \ref{theo:treeres}. However, a simple calculation shows that in $\mathcal{G}^\text{star}_3$,
\[
r_{2,3} = 2\left(n + m\right) - \frac{2nm}{n+m},
\]
which only matches \refe{eqn:treeres} for $n = m = 1$. Thus in more general cases of connection subgraphs like $\mathcal{G}^\text{tree}_{n,m}$ but with arbitrary weights on every edge, the resistance between the leaves does \emph{not} depend only on the equivalent resistance of each path.

Theorems \ref{theo:pathres}, \ref{theo:cycleres} and \ref{theo:treeres} by no means characterise all the possible connection subgraphs in a directed graph. Other connection subgraphs include multiple paths from $k$ to $j$ (some of which could coincide over part of their length), multiple paths from $k$ to $j$ and multiple paths from $j$ to $k$ (again, some of which could partially coincide), multiple indirect connections of the type analysed in Theorem \ref{theo:treeres} (which could partially coincide) and a combination of indirect and direct (i.e. path) connections. Further analysis is needed to completely describe how to compute resistances in these situations.

\section{Conclusions}\label{sec:conc}
The results of Lemma \ref{lem:edgeres} and Theorems \ref{theo:pathres} and \ref{theo:cycleres} demonstrate that in some situations our definition of effective resistance for directed graphs behaves as an intuitive extension of effective resistance in undirected graphs. In contrast, Theorem \ref{theo:treeres} demonstrates a fundamental difference between effective resistance in directed and undirected graphs that arises from the fundamentally different connections that are possible only in directed graphs. Nevertheless, the results presented above show that our notion of effective resistance for directed graphs provides an approach that can relate the local structure of a directed graph to its global properties. The familiar properties of effective resistance allows for a firm analysis of directed graphs that behave similarly to undirected graphs, while the unfamiliar properties can provide insight for the design of directed networks which contain essential differences as compared to undirected networks.

\FloatBarrier


\appendices

\section{Proof of Lemma \ref{lem:treeresnl1}}\label{apndx:treeres}
\begin{IEEEproof}
As stated in the lemma, we will assume that $n$ and $\ell$ are positive integers throughout this proof. Let $N_{n,\ell}$ be the number of nodes in $\mathcal{G}^\text{tree}_{n,\ell}$. The branch of length $n$ contains $n$ nodes (excluding the root), while the other branch contains $\ell$ nodes (excluding the root). Therefore, we have $N_{n,\ell} = n+\ell+1$. Let us label the nodes in $\mathcal{G}^\text{tree}_{n,\ell}$ from $1$ to $n+1$ along the branch of length $n$, in reverse order of the edge directions and starting with the root (thus the leaf of this branch is node $n+1$). Then let us label the nodes in the branch of length $\ell$ from $n+2$ to $N_{n,\ell} = n+\ell+1$ in reverse order of the edge directions. Thus the second leaf is node $N_{n,\ell}$.

In the following, we will denote the adjacency matrix of $\mathcal{G}^\text{tree}_{n,\ell}$ by $A_{n,\ell}$, its matrix of node out-degrees by $D_{n,\ell}$ and its Laplacian matrix by $L_{n,\ell}$. Furthermore, we will let $Q_{n,\ell}$ be a $\left(N_{n,\ell}-1\right)\times N_{n,\ell}$ matrix that satisfies \refe{eqn:propq} and $\overline{L}_{n,\ell}$ and $\Sigma_{n,\ell}$ be the corresponding matrices from \refe{eqn:lbar} and \refe{eqn:lyap} using $L_{n,\ell}$ and $Q_{n,\ell}$. Finally, $X_{n,\ell}$ will be the matrix from \refe{eqn:xdef}, computed using $\Sigma_{n,\ell}$ and $Q_{n,\ell}$. Then, by Lemma \ref{lem:Xres}, the entries of $X_{n,\ell}$ are related to the resistances in $\mathcal{G}^\text{tree}_{n,\ell}$ by \refe{eqn:rtox}.

As in the proof of Lemma \ref{lem:treeresn1}, let $A_{N_p}^\text{path}$, $D_{N_p}^\text{path}$ and $L_{N_p}^\text{path}$ denote the adjacency matrix, matrix of out-degrees and Laplacian matrix of a directed path containing $N_p$ nodes and unit weights on every edge. Let the nodes in this path be labelled from $1$ to $N_p$ in the order in which they appear, starting with the root. Then we can write $A_{n,\ell}$, $D_{n,\ell}$ and $L_{n,\ell}$ in terms of $A_{N_p}^\text{path}$, $D_{N_p}^\text{path}$ and $L_{N_p}^\text{path}$ as follows: $A_{n,\ell} = \begin{bmatrix}A_{n+1}^\text{path} & \Mzero\\ \mathbf{e}_{\ell}^{(1)}\mathbf{e}_{n+1}^{(1) T} & A_{\ell}^\text{path}\end{bmatrix}$, $D_{n,\ell} = \begin{bmatrix}D_{n+1}^\text{path} & \Mzero\\ \Mzero & D_{\ell}^\text{path} + \mathbf{e}_{\ell}^{(1)}\mathbf{e}_{\ell}^{(1) T}\end{bmatrix}$ and $L_{n,\ell} = \begin{bmatrix}L_{n+1}^\text{path} & \Mzero\\ -\mathbf{e}_{\ell}^{(1)}\mathbf{e}_{n+1}^{(1) T} & L_{\ell}^\text{path} + \mathbf{e}_{\ell}^{(1)}\mathbf{e}_{\ell}^{(1) T}\end{bmatrix}$.

Using these expressions as well as \refe{eqn:oneL} and \refe{eqn:eiL}, we can observe that
\begin{equation}\label{eqn:oneLnl}
\mathbf{1}_{N_{n,\ell}}^T L_{n,\ell} = -2\mathbf{e}_{N_{n,\ell}}^{(1) T} + \mathbf{e}_{N_{n,\ell}}^{(n+1) T} + \mathbf{e}_{N_{n,\ell}}^{(N_{n,\ell}) T} \text{, and}
\end{equation}
\begin{equation}\label{eqn:eiLnl}
\mathbf{e}_{N_{n,\ell}}^{(i)T} L_{n,\ell} = \begin{cases}\mathbf{e}_{N_{n,\ell}}^{(i) T} \!- \mathbf{e}_{N_{n,\ell}}^{(i-1) T} &\text{if } 1 < i \leq N_{n,\ell}, \, i \neq n+2,\\\mathbf{e}_{N_{n,\ell}}^{(n+2) T} \!- \mathbf{e}_{N_{n,\ell}}^{(1) T} &\text{if } i = n+2,\\\mathbf{0}^T &\text{if } i = 1.\end{cases}
\end{equation}

Let us now consider $\mathcal{G}^\text{tree}_{n,\ell+1}$. By our labeling convention, the resistance between the two leaves of $\mathcal{G}^\text{tree}_{n,\ell+1}$ is given by $r(n,\ell+1) = r_{n+1,n+\ell+2}$. Now, we can write the adjacency matrix of $\mathcal{G}^\text{tree}_{n,\ell+1}$ in terms of $A_{n,\ell}$ as $A_{n,\ell+1} = \begin{bmatrix}A_{n,\ell} & \mathbf{0}\\\mathbf{e}_{N_{n,\ell}}^{(N_{n,\ell}) T} & 0\end{bmatrix}$. In a similar fashion, we can write the matrix of node out-degrees for $\mathcal{G}^\text{tree}_{n,\ell+1}$ as $D_{n,\ell+1} = \begin{bmatrix}D_{n,\ell} & \mathbf{0}\\\mathbf{0}^T & 1\end{bmatrix}$, and the Laplacian matrix as $L_{n,\ell+1} = \begin{bmatrix}L_{n,\ell} & \mathbf{0} \\ -\mathbf{e}_{N_{n,\ell}}^{(N_{n,\ell}) T} & 1\end{bmatrix}$.

Now, let $Q_{n,\ell+1} = \begin{bmatrix}Q_{n,\ell} & \mathbf{0}\\\alpha\mathbf{1}_{N_{n,\ell}}^T & -\beta\end{bmatrix}$, where $\alpha = \frac{1}{\sqrt{N_{n,\ell}(N_{n,\ell}+1)}} = \frac{1}{\sqrt{(n+\ell+1)(n+\ell+2)}}$ and $\beta = \sqrt{\frac{N_{n,\ell}}{N_{n,\ell}+1}} = \sqrt{\frac{n+\ell+1}{n+\ell+2}}$. Then $Q_{n,\ell+1}$ satisfies \refe{eqn:propq}. We can therefore use \refe{eqn:lbar}, \refe{eqn:oneLnl} and the facts that $L_{n,\ell} \mathbf{1}_{N_{n,\ell}} = \mathbf{0}_{N_{n,\ell}}$ and $\beta\left(\alpha + \beta\right) = 1$ to compute $\overline{L}_{n,\ell}$ as
\begin{multline*}
\overline{L}_{n,\ell} \!\!=\!\!\! \left[\begin{matrix}\overline{L}_{n,\ell} \\  \!-2\alpha\mathbf{e}_{N_{n,\ell}}^{(1) T}Q_{n,\ell}^T \!+\! \alpha\mathbf{e}_{N_{n,\ell}}^{(n+1) T}Q_{n,\ell}^T \!+\! \left(\alpha \!+\! \beta\right)\mathbf{e}_{N_{n,\ell}}^{(N_{n,\ell}) T}Q_{n,\ell}^T \end{matrix}\right.\\
\left.\begin{matrix} \mathbf{0}\\  1\end{matrix}\right].
\end{multline*}

In order to compute resistances in $\mathcal{G}^\text{tree}_{n,\ell+1}$, we must find the matrix $\Sigma_{n,\ell+1}$ which solves \refe{eqn:lyap}. Since we have partitioned $\overline{L}_{n,\ell+1}$ into a $2\times 2$ block matrix, we will do the same for $\Sigma_{n,\ell+1}$. Let $\Sigma_{n,\ell+1} = \begin{bmatrix} S & \mathbf{t}\\ \mathbf{t}^T & u\end{bmatrix}$, where $S \in \mathbb{R}^{(N_{n,\ell}-1)\times (N_{n,\ell}-1)}$ is a symmetric matrix, $\mathbf{t} \in \mathbb{R}^{N_{n,\ell}-1}$ and $u \in \mathbb{R}$. Then multiplying out the matrices in \refe{eqn:lyap} and equating blocks in this matrix equation gives us
\begin{align}
\overline{L}_{n,\ell} S + S\overline{L}_{n,\ell}^T						&= I_{N_{n,\ell}-1}, \label{eqn:Tnl1lyapS} \\
\overline{L}_{n,\ell} \mathbf{t} + \mathbf{t} - 2\alpha S Q_{n,\ell} \mathbf{e}_{N_{n,\ell}}^{(1)} + \alpha S Q_{n,\ell}\mathbf{e}_{N_{n,\ell}}^{(n+1)} \hspace{-0.2cm}& \nonumber\\
{} + \left(\alpha + \beta\right) S Q_{n,\ell}\mathbf{e}_{N_{n,\ell}}^{(N_{n,\ell})}	&= \mathbf{0} \text{, and} \label{eqn:Tnl1lyapt} \\
2u - 4\alpha\mathbf{e}_{N_{n,\ell}}^{(1) T}Q_{n,\ell}^T\mathbf{t} + 2\alpha \mathbf{e}_{N_{n,\ell}}^{(n+1) T}Q_{n,\ell}^T\mathbf{t} \hspace{0.6cm}& \nonumber\\
{} + 2\left(\alpha + \beta\right) \mathbf{e}_{N_{n,\ell}}^{(N_{n,\ell}) T}Q_{n,\ell}^T\mathbf{t}			&= 1. \label{eqn:Tnl1lyapu}
\end{align}

From \refe{eqn:Tnl1lyapS}, it is clear that $S = \Sigma_{n,\ell}$. In addition, we can rewrite \refe{eqn:Tnl1lyapu} as
\begin{multline}\label{eqn:Tnl1uform}
u = \frac{1}{2} + 2\alpha\mathbf{e}_{N_{n,\ell}}^{(1) T}Q_{n,\ell}^T\mathbf{t} - \alpha \mathbf{e}_{N_{n,\ell}}^{(n+1) T}Q_{n,\ell}^T\mathbf{t} \\
{} - \left(\alpha + \beta\right) \mathbf{e}_{N_{n,\ell}}^{(N_{n,\ell}) T}Q_{n,\ell}^T\mathbf{t}.
\end{multline}
Thus in order to find a complete solution for $\Sigma_{n,\ell+1}$, we must solve \refe{eqn:Tnl1lyapt} for $\mathbf{t}$. However, resistances are computed from the entries of $X_{n,\ell+1}$, which, if we let $\mathbf{v} \mathrel{\mathop :}= Q_{n,\ell}^T \mathbf{t} = \left[v_i\right]$ and use \refe{eqn:xdef}, can be written as
\begin{multline*}
X_{n,\ell+1} \!=\! \left[\begin{matrix} X_{n,\ell} + 2\alpha \mathbf{v}\mathbf{1}_{N_{n,\ell}}^T \!\!+ 2\alpha\mathbf{1}_{N_{n,\ell}} \mathbf{v}^T \!\!+ 2\alpha^2 u \mathbf{1}_{N_{n,\ell}}\mathbf{1}_{N_{n,\ell}}^T \\ -2\beta \mathbf{v}^T - 2\alpha\beta u \mathbf{1}_{N_{n,\ell}}^T\end{matrix}\right.\\
\left.\begin{matrix} -2\beta \mathbf{v} - 2\alpha\beta u \mathbf{1}_{N_{n,\ell}}\\  2\beta^2 u \end{matrix}\right].
\end{multline*}
Hence, in order to compute resistances in $\mathcal{G}^\text{tree}_{n,\ell+1}$, we need only compute $\mathbf{v}$, not $\mathbf{t}$. We should also note that as $X_{n,\ell+1}$ does not depend on our choice of $Q_{n,\ell+1}$ (by Lemma \ref{P1:lem:indofq} in the companion paper \cite{Young13I}), neither does $\mathbf{v}$. In fact, we can write \refe{eqn:Tnl1uform} as $u = \frac{1}{2} + 2\alpha v_{1} - \alpha v_{n+1} - \left(\alpha + \beta\right) v_{N_{n,\ell}}$, and the resistance we seek as
\begin{multline}
r(n,\ell+1) = x_{n,\ell \, n+1,n+1} + \left(\alpha +\beta\right)^2 + 4\alpha\left(\alpha +\beta\right)^2 v_1 \\
{} + 2\left(\alpha + \beta\right)\left[2 - \alpha\left(\alpha + \beta\right)\right] v_{n+1} - 2(\alpha + \beta)^3 v_{N_{n,\ell}}. \label{eqn:rnl1fromv}
\end{multline}
Thus we only need to find $v_1$, $v_{n+1}$ and $v_{N_{n,\ell}}$ in order to compute $r(n,\ell+1)$.

Now, $v_i = \mathbf{e}_{N_{n,\ell}}^{(i) T}\mathbf{v} = \mathbf{e}_{N_{n,\ell}}^{(i) T} Q_{n,\ell}^T \mathbf{t}$. We will therefore proceed by left-multiplying \refe{eqn:Tnl1lyapt} by $\mathbf{e}_{N_{n,\ell}}^{(i) T} Q_{n,\ell}^T$. Using the fact that $S = \Sigma_{n,\ell}$, we obtain
\begin{multline}\label{eqn:Tnl1lyapv}
\mathbf{e}_{N_{n,\ell}}^{(i) T} Q_{n,\ell}^T Q_{n,\ell} L_{n,\ell} \mathbf{v} + v_i - \alpha \mathbf{e}_{N_{n,\ell}}^{(i) T} X_{n,\ell} \mathbf{e}_{N_{n,\ell}}^{(1)} \\
{}+ \frac{\alpha}{2} \mathbf{e}_{N_{n,\ell}}^{(i) T} X_{n,\ell} \mathbf{e}_{N_{n,\ell}}^{(n+1)} + \frac{\alpha + \beta}{2} \mathbf{e}_{N_{n,\ell}}^{(i) T} X_{n,\ell} \mathbf{e}_{N_{n,\ell}}^{(N_{n,\ell})} = 0.
\end{multline}
But $\mathbf{e}_{N_{n,\ell}}^{(i) T} Q_{n,\ell}^T Q_{n,\ell} = \mathbf{e}_{N_{n,\ell}}^{(i) T} \left(I_{N_{n,\ell}} - \frac{1}{N_{n,\ell}}\mathbf{1}_{N_{n,\ell}}\mathbf{1}_{N_{n,\ell}}^T\right) = \mathbf{e}_{N_{n,\ell}}^{(i) T} - \frac{1}{N_{n,\ell}}\mathbf{1}_{N_{n,\ell}}^T$ by \refe{eqn:propq}, and so by using \refe{eqn:oneLnl} and \refe{eqn:eiLnl}, we find
\begin{multline*}
\mathbf{e}_{N_{n,\ell}}^{(i) T} Q_{n,\ell}^T Q_{n,\ell}L_{n,\ell}\mathbf{v} = \frac{2}{N_{n,\ell}}v_1 - \frac{1}{N_{n,\ell}}v_{n+1} - \frac{1}{N_{n,\ell}}v_{N_{n,\ell}} \\
{} + \begin{cases}v_i - v_{i-1} & \text{if } 1 < i \leq N_{n,\ell}, \; i \neq n+2,\\ v_{n+2} - v_{1} & \text{if } i = n+2,\\ 0 & \text{if } i = 1.\end{cases}
\end{multline*}
Furthermore, we observe that $\mathbf{e}_{N_{n,\ell}}^{(i) T} X_{n,\ell} \mathbf{e}_{N_{n,\ell}}^{(1)} = x_{n,\ell \, i,1}$, $\mathbf{e}_{N_{n,\ell}}^{(i) T} X_{n,\ell} \mathbf{e}_{N_{n,\ell}}^{(n+1)} = x_{n,\ell \, i,n+1}$, and $\mathbf{e}_{N_{n,\ell}}^{(i) T} X_{n,\ell} \mathbf{e}_{N_{n,\ell}}^{(N_{n,\ell})} = x_{n,\ell \, i,N_{n,\ell}}$.

Substituting these expressions into \refe{eqn:Tnl1lyapv} gives us
\begin{align}
v_i		&= \frac{1}{2}v_{i-1} - \frac{1}{N_{n,\ell}}v_1 + \frac{1}{2N_{n,\ell}}v_{n+1} + \frac{1}{2N_{n,\ell}}v_{N_{n,\ell}}  \nonumber\\
		& \hspace{0.0cm} {} + \frac{\alpha}{2}x_{n,\ell\, i,1} - \frac{\alpha}{4}x_{n,\ell\, i,n+1} - \frac{\alpha+\beta}{4}x_{n,\ell\, i,N_{n,\ell}}  \label{eqn:recurrnlvi}\\
		& \hspace{4cm }\text{ if } 1< i \leq N_{n,\ell}, \; i \neq n+2, \nonumber\\
v_{n+2}	&= \frac{1}{2}v_{1} - \frac{1}{N_{n,\ell}}v_1 + \frac{1}{2N_{n,\ell}}v_{n+1} + \frac{1}{2N_{n,\ell}}v_{N_{n,\ell}} \nonumber\\
		& \hspace{1.6cm} {} + \frac{\alpha}{2}x_{n,\ell\, n+2,1} - \frac{\alpha}{4}x_{n,\ell\, n+2,n+1} \nonumber\\
		& \hspace{3.2cm} {} - \frac{\alpha+\beta}{4}x_{n,\ell\, n+2,N_{n,\ell}} \text{, and} \label{eqn:recurrnlvn2} \\
v_1		&= \frac{1}{N_{n,\ell}+2}v_{n+1} + \frac{1}{N_{n,\ell} + 2}v_{N_{n,\ell}} + \frac{\alpha N_{n,\ell}}{N_{n,\ell} + 2}x_{n,\ell\, 1,1} \nonumber\\
		& \hspace{-0.5cm} {} - \frac{\alpha N_{n,\ell}}{2\left(N_{n,\ell} + 2\right)}x_{n,\ell\, 1,n+1} - \frac{\left(\alpha + \beta\right)N_{n,\ell}}{2\left(N_{n,\ell} + 2\right)}x_{n,\ell\, 1,N_{n,\ell}}.\label{eqn:recurrnlv1}
\end{align}

We can now recursively apply \refe{eqn:recurrnlvi} $n$ times, starting with $i = n+1$, and simplify using \refe{eqn:sumtwos} to find
\begin{multline}
\frac{N_{n,\ell} \!-\! 1 \!+\! 2^{-n}}{N_{n,\ell}}v_{n+1} \!=\! \frac{\left[-2 \!+\! \left(N_{n,\ell} \!+\! 2\right)2^{-n}\right]}{N_{n,\ell}}v_1 + \frac{1 \!-\! 2^{-n}}{N_{n,\ell}}v_{N_{n,\ell}} \\
{} + \alpha\sum_{k = 1}^{n}{x_{n,\ell\, n+2-k,1} 2^{-k}} - \frac{\alpha}{2}\sum_{k = 1}^{n}{x_{n,\ell\, n+2-k,n+1} 2^{-k}} \\
	{} - \frac{\alpha + \beta}{2}\sum_{k = 1}^{n}{x_{n,\ell\, n+2-k,N_{n,\ell}} 2^{-k}}. \label{eqn:recurrnlvn1}
\end{multline}

Similarly, we can recursively apply \refe{eqn:recurrnlvi} $\ell-1$ times, starting with $i = N_{\ell,n} = n+\ell+1$, substitute in \refe{eqn:recurrnlvn2} and simplify using \refe{eqn:sumtwos} to find
\begin{multline}
\frac{N_{n,\ell} \!-\! 1 \!+\! 2^{-\ell}}{N_{n,\ell}}v_{N_{n,\ell}} \!=\! \frac{\left[-2 \!+\! \left(N_{n,\ell} \!+\! 2\right)2^{-\ell}\right]}{N_{n,\ell}}v_1 + \frac{1 \!-\! 2^{-\ell}}{N_{n,\ell}}v_{n+1} \\
{} + \alpha\sum_{k = 1}^{\ell}{x_{n,\ell\, N_{n,\ell}+1-k,1} 2^{-k}}  - \frac{\alpha}{2}\sum_{k = 1}^{\ell}{x_{n,\ell\, N_{n,\ell}+1-k,n+1} 2^{-k}} \\
	{} - \frac{\alpha + \beta}{2}\sum_{k = 1}^{\ell}{x_{n,\ell\, N_{n,\ell}+1-k,N_{n,\ell}} 2^{-k}}. \label{eqn:recurrnlvN}
\end{multline}
Note that \refe{eqn:recurrnlvN} reduces to \refe{eqn:recurrnlvn2} when $\ell = 1$.

But now \refe{eqn:recurrnlv1}, \refe{eqn:recurrnlvn1} and \refe{eqn:recurrnlvN} form a set of three of simultaneous linear equations in $v_1$, $v_{n+1}$ and $v_{N_{n,\ell}}$. Substituting their solution into \refe{eqn:rnl1fromv} and then multiplying by $N_{n,\ell}$ (and using the definitions of $\alpha$ and $\beta$) gives us
\begin{multline}\label{eqn:rnl1fromx}
N_{n,\ell}\, r(n,\ell+1) = N_{n,\ell} + 1 + \left(2^{2-n} - 2^{1-\ell}\right)x_{n,\ell\, 1,1} \\
{} + \left(2^{-\ell} \!-\! 2^{1-n}\right)x_{n,\ell\, 1,n+1} + \left(N_{n,\ell} \!+\! 1\right)\left(2^{-\ell} \!-\! 2^{1-n}\right)x_{n,\ell\, 1,N_{n,\ell}}\\
	{}  + N_{n,\ell}x_{n,\ell\, n+1,n+1} + 4\sum_{k = 1}^{n}{x_{n,\ell\, n+2-k,1} 2^{-k}} \\
	{} - 2\sum_{k = 1}^{\ell}{x_{n,\ell\, N_{n,\ell}+1-k,1} 2^{-k}} - 2\sum_{k = 1}^{n}{x_{n,\ell\, n+2-k,n+1} 2^{-k}} \\
	{} \!+\! \sum_{k = 1}^{\ell}{\!x_{n,\ell\, N_{n,\ell}+1-k,n+1} 2^{-k}} \!- \left(2N_{n,\ell} \!+\! 2\right)\!\!\sum_{k = 1}^{n}{\!x_{n,\ell\, n+2-k,N_{n,\ell}} 2^{-k}} \\
	{}+ \left(N_{n,\ell} + 1\right)\sum_{k = 1}^{\ell}{x_{n,\ell\, N_{n,\ell}+1-k,N_{n,\ell}} 2^{-k}}.
\end{multline}

Now, by \refe{eqn:dirres}, we can write $x_{n,\ell\, k,j} = \frac{1}{2}x_{n,\ell\, k,k} + \frac{1}{2}x_{n,\ell\, j,j} - \frac{1}{2}r_{k,j}$. Furthermore, by Theorem~\ref{theo:pathres} we know that
\begin{multline}\label{eqn:rpathintree}
r_{k,j} = 2\left|k - j\right| \text{ if } 1\leq k,j \leq n+1, \\
\text{or } n+2 \leq k,j \leq N_{n,\ell} \text{, and}
\end{multline}
\begin{equation}\label{eqn:rpath2intree}
r_{1,j} = 2\left(j - n - 1\right) \text{ if } n+2 \leq j \leq N_{n,\ell}.
\end{equation}
Finally, by the definition of $r(\cdot,\cdot)$, we can say that
\begin{multline}\label{eqn:rcrossintree}
r_{k,j} = r(k-1,j-n-1) \text{ if } 1 < k \leq n+1 \\
\text{ and } n+2 \leq j \leq N_{n,\ell}.
\end{multline}
Therefore, we can substitute for each non-diagonal $x_{n,\ell\, k,j}$ term in \refe{eqn:rnl1fromx} and use \refe{eqn:sumtwos} and \refe{eqn:suminttwos}, along with the fact that $N_{n,\ell} = n + \ell + 1$ to find
\begin{multline*}
\left(n\!+\!\ell\!+\!1\right) r(n,\ell+1) \!=\! -4n + 2\ell + 4 + \left(\ell^2 + n\ell + 2\ell - 3\right)\!2^{1-n} \\
	{} + \left(\ell + 4\right)2^{-\ell} + \left(n \!+\! \ell \!+\! 2\right)\sum_{k=1}^{n}{r(n+1-k,\ell)2^{-k}} \\
	{} - \frac{1}{2}\sum_{k=1}^{\ell}{r(n,\ell+1-k)2^{-k}} + \left(n \!+\! \ell \!+\! \frac{1}{2}\right)x_{n,\ell\, n+1,n+1}\\
	{} + \left[\left(n+\ell+1\right)\left(2^{-1-\ell} - 2^{-n}\right) + 1\right]x_{n,\ell\, 1,1} \\
	{} - \frac{n+\ell+2}{2}x_{n,\ell\, n+\ell+1,n+\ell+1} \\
	{}  - \left(n + \ell + 1\right)\!\sum_{k=1}^{n}{x_{n,\ell \, n+2-k,n+2-k}2^{-k}} \\
	{} + \frac{n + \ell + 1}{2}\!\sum_{k=1}^{\ell}{x_{n,\ell \, n+\ell+2-k,n+\ell+2-k}2^{-k}},
\end{multline*}
or, by changing indices inside the sums,
\begin{multline}\label{eqn:rnl1fromdiagx}
\left(n\!+\!\ell\!+\!1\right) r(n,\ell+1) \!=\! -4n + 2\ell + 4 + \left(\ell^2 + n\ell + 2\ell - 3\right)\!2^{1-n} \\
	{} + \left(\ell + 4\right)2^{-\ell} + \frac{n + \ell + 2}{2}\sum_{k=1}^{n}{r(k,\ell)2^{k-n}} - \frac{1}{4}\sum_{k=1}^{\ell}{r(n,k)2^{k-\ell}} \\
	{} + \left(n + \ell + \frac{1}{2}\right)x_{n,\ell\, n+1,n+1} \\
	{} + \left[\left(n+\ell+1\right)\left(2^{-1-\ell} - 2^{-n}\right) + 1\right]x_{n,\ell\, 1,1} \\
	{} \!-\! \frac{n\!+\!\ell\!+\!2}{2}x_{n,\ell\, n+\ell+1,n+\ell+1} \!-\! \frac{n \!+\! \ell \!+\! 1}{2}\sum_{k=1}^{n}{\!x_{n,\ell \, k+1,k+1}2^{k-n}}\\
	{} + \frac{n + \ell + 1}{4}\sum_{k=1}^{\ell}{x_{n,\ell \, n+1+k,n+1+k}2^{k-\ell}}.
\end{multline}

Now, by \refe{eqn:rtox} from Lemma \ref{lem:Xres}, we know that
\begin{equation}\label{eqn:xnldiagform}
x_{n,\ell\, i,i} = \frac{1}{N_{n,\ell}}\sum_{k = 1}^{N_{n,\ell}}{r_{i,k}} - \frac{1}{N_{n,\ell}^2}\sum_{k = 1}^{N_{n,\ell}-1}{\sum_{j = k+1}^{N_{n,\ell}}{r_{k,j}}}.
\end{equation}
Using \refe{eqn:rpathintree}, \refe{eqn:rpath2intree} and \refe{eqn:rcrossintree} and then \refe{eqn:sumint}, we can write the first sum in \refe{eqn:xnldiagform} as
\begin{equation}\label{eqn:nlsumresrow}
\sum_{k = 1}^{n+\ell+1}{\!\!\!r_{i,k}} = \begin{cases}n^2 + \ell^2 + n + \ell & \text{if } i = 1,\\ n^2 + (3-2i)n + 2(i-1)^2 + \displaystyle\sum_{k=1}^{\ell}{r(i-1,k)}\hspace{-4cm}&\\ & \text{if } 1 < i \leq n+1,\\ 2n^2 + \ell^2 + 2n\ell + (4-4i)n + (3-2i)\ell \hspace{-4cm}& \\ {} + 2(i-1)^2 + \displaystyle\sum_{k = 1}^{n}{r(k,i-n-1)} \hspace{-3cm}&\\& \text{if } n+2 \leq i \leq n+\ell+1.\end{cases}
\end{equation}

We can note that the double sum in \refe{eqn:xnldiagform} is independent of $i$. Let
\[
f \mathrel{\mathop :}= \sum_{k = 1}^{N_{n,\ell}-1}{\sum_{j = k+1}^{N_{n,\ell}}{r_{k,j}}}.
\]
Then, substituting \refe{eqn:xnldiagform} and \refe{eqn:nlsumresrow} into \refe{eqn:rnl1fromdiagx}, changing indices and using the results of Lemma \ref{lem:finser} produces
\begin{multline}\label{eqn:Nrnl1}
\left(n+\ell+1\right) r(n,\ell+1) = \frac{-3n^2 + 3\ell^2 - 2n\ell - n + 5\ell + 2}{2(n+\ell+1)} \\
	{} + \left(\ell^2 + 2n\ell + 2n + 3\ell\right)2^{-n} + \left(n^2 + n + 2\right)2^{-1-\ell} \\
	{} + \frac{1}{4}\sum_{k=1}^{\ell}{\left(4 - \frac{2}{n+\ell+1} - 2^{k-\ell}\right)r(n,k)}\\
	{} - \frac{n+\ell+2}{2}\sum_{k=1}^{n}{\left(\frac{1}{n+\ell+1} - 2^{k-n}\right)r(k,\ell)} \\
	{} - \frac{1}{4}\sum_{k=1}^{n}{\sum_{j=1}^{\ell}{\left(2^{1+k-n} - 2^{j-\ell}\right)r(k,j)}}.
\end{multline}
Finally, dividing \refe{eqn:Nrnl1} through by $n + \ell + 1$ produces our desired result.
\end{IEEEproof}

\section{Finite series}\label{apndx:finser}
The following series are either well-known or special cases of well-known series. The first two and the general cases of the third and fourth usually appear in any introductory mathematical text that covers series (e.g. section 4.2 of \cite{Riley06}). The fifth is slightly more obscure.
\begin{lemma}\label{lem:finser}
For integer values of $n > 0$,
\begin{enumerate}[(i) ]
\item \begin{bulletequation}\label{eqn:sumint} \sum_{k = 1}^{n}{k} = \frac{1}{2}n(n+1), \end{bulletequation}

\item \begin{bulletequation}\label{eqn:sumintsq} \sum_{k = 1}^{n}{k^2} = \frac{1}{6}n(n+1)(2n+1), \end{bulletequation}

\item \begin{bulletequation}\label{eqn:sumtwos} \sum_{k = 1}^{n}{2^{-k}} = 1 - 2^{-n}, \end{bulletequation}

\item \begin{bulletequation}\label{eqn:suminttwos} \sum_{k = 1}^{n}{k2^{-k}} = 2 - (n+2)2^{-n} \text{, and} \end{bulletequation}

\item \begin{bulletequation}\label{eqn:sumintsqtwos} \sum_{k = 1}^{n}{k^2 2^{-k}} = 6 - \left(n^2 + 4n + 6\right)2^{-n}. \end{bulletequation}

\end{enumerate}
\end{lemma}

\begin{IEEEproof}
Equations \refe{eqn:sumint} and \refe{eqn:sumintsq} are special cases of (6.2.1) in \cite{Hansen75}, while \refe{eqn:sumtwos}, \refe{eqn:suminttwos} and \refe{eqn:sumintsqtwos} are special cases of (6.9.1) in \cite{Hansen75}. All are easily proved using induction.
\end{IEEEproof}

\subsection{Finite series of binomial coefficients}\label{subsec:binsums}
Although there are many interpretations and uses of binomial coefficients, we will simply assume two basic facts about them, namely \emph{Pascal's rule};
\begin{equation}\label{eqn:pascal}
\binom{n}{k} = \binom{n-1}{k} + \binom{n-1}{k-1}, \; 1 \leq k \leq n-1,
\end{equation}
and the \emph{binomial formula};
\begin{equation}\label{eqn:binform}
\left(x + y\right)^n = \sum_{i = 0}^n{\binom{n}{i} x^i y^{n-i}}, \; n \geq 0.
\end{equation}
Pascal's rule follows easily from \refe{eqn:bincoef} while the binomial formula can be inductively proved using Pascal's rule. Equations \refe{eqn:pascal} and \refe{eqn:binform} can also be found in standard introductory mathematics texts, such as sections 1.5--1.6 in \cite{Riley06}.

We can use Pascal's rule to derive some identities involving binomial coefficients. These identities include the two in the following lemma.

\begin{lemma}\label{lem:binsumnmk}
For integer values of $n, m$ and $k$, with $n > 0$, $m \geq 0$ and $0 \leq k \leq m$,
\begin{enumerate}[(i) ]
\item \begin{bulletequation}\label{eqn:binsumnm}\sum_{i=1}^{n}{\binom{m+i}{m+1}} = \binom{n+m+1}{m+2} \text{, and}\end{bulletequation}

\item \begin{bulletequation}\label{eqn:binsumnmk}\sum_{i=1}^{n}{\binom{m+i}{k+i}} = \binom{n+m+1}{n+k} - \binom{m+1}{k}.\end{bulletequation}
\end{enumerate}
\end{lemma}

\begin{IEEEproof}
Both results can be easily proven using mathematical induction and Pascal's rule.
\end{IEEEproof}

A special case of the binomial formula can be found by substituting $y = 1$ into \refe{eqn:binform}, which gives
\begin{equation}\label{eqn:binformx}
\left(1 + x\right)^n = \sum_{i = 0}^n{\binom{n}{i} x^i}, \; n \geq 0.
\end{equation}
Differentiating this expression with respect to $x$ gives us
\begin{equation}\label{eqn:binformder}
n\left(1 + x\right)^{n-1} = \sum_{i = 0}^n{i\binom{n}{i} x^{i-1}}, \; n \geq 1.
\end{equation}

In the following results, we will make use of a few ``well-known'' series of binomial coefficients (for example, the first two can be found in Chapter 3 of \cite{Spiegel09} and all can be solved by Mathematica). Since they are not as standard as the basic facts stated above, we will include a brief proof of them for the sake of completeness.

\begin{lemma}[Standard sums of binomial coefficients]\label{lem:standardsums}
For integer values of $n$,
\begin{enumerate}[(i) ]
\item \begin{bulletequation}\label{eqn:binevensum}\sum_{i = 0}^{\Floor{\frac{n}{2}}}{\binom{n}{2i}} = 2^{n-1}, \; n > 0,\end{bulletequation}

\item \begin{bulletequation}\label{eqn:binoddsum}\sum_{i = 0}^{\Floor{\frac{n-1}{2}}}{\binom{n}{2i+1}} = 2^{n-1}, \; n > 0,\end{bulletequation}

\item \begin{bulletequation}\label{eqn:binevensumi}\sum_{i = 0}^{\Floor{\frac{n}{2}}}{2i\binom{n}{2i}} = n2^{n-2}, \; n > 1\text{, and} \end{bulletequation}

\item \begin{bulletequation}\label{eqn:binoddsumi}\sum_{i = 0}^{\Floor{\frac{n-1}{2}}}{\left(2i+1\right)\binom{n}{2i+1}} = n2^{n-2}, \; n > 1.\end{bulletequation}

\end{enumerate}
\end{lemma}

\begin{IEEEproof}
Substituting $x = \pm 1$ into \refe{eqn:binformx} gives us $\sum_{i = 0}^n{\binom{n}{i}} = 2^n$ and $\sum_{i = 0}^n{(-1)^i \binom{n}{i}} = 0$ for any $n > 0$. Equations \refe{eqn:binevensum} and \refe{eqn:binoddsum} can be found by taking the sum and difference of these two expressions and dividing by $2$. Similarly, substituting $x = \pm 1$ into \refe{eqn:binformder} gives us $\sum_{i = 0}^n{i\binom{n}{i}} = n2^{n-1}$ and $\sum_{i = 0}^n{(-1)^{i-1}i \binom{n}{i}} = 0$ for any $n > 1$. Equations \refe{eqn:binevensumi} and \refe{eqn:binoddsumi} can be found by taking the sum and difference of these two expressions and dividing by $2$.
\end{IEEEproof}

We can now use the results from Lemma \ref{lem:standardsums} to derive some more specialised series. These are summarised in the following lemma. As a point of notation, we will assume that any sum not containing any terms (such as $\sum_{i=0}^{-1}a_i$) is equal to zero.

\begin{lemma}[Specialised sums of binomial coefficients] \label{lem:specialsums}
For integer values of $p \geq 0$,
\begin{enumerate}[(i) ]

\item \begin{bulletequation}\label{eqn:ievenevensum4}\sum_{i = 1}^{p+1}{i\binom{2p+4}{2i+2}} = p2^{2p+2} + 1,\end{bulletequation}

\item \begin{bulletequation}\label{eqn:ievenoddsum}\sum_{i = 1}^{p}{i\binom{2p+2}{2i+1}} = p2^{2p},\end{bulletequation}

\item \begin{bulletequation}\label{eqn:ievenoddsum4}\sum_{i = 1}^{p+1}{i\binom{2p+4}{2i+1}} = \left(p+1\right)2^{2p+2},\end{bulletequation}

\item \begin{bulletequation}\label{eqn:ioddevensum}\sum_{i = 1}^{p}{i\binom{2p+3}{2i+2}} = \left(2p-1\right)2^{2p} + 1,\end{bulletequation}

\item \begin{bulletequation}\label{eqn:ioddoddsum}\sum_{i = 1}^{p}{i\binom{2p+3}{2i+1}} = \left(2p+1\right)2^{2p} - p - 1,\end{bulletequation}

\item \begin{bulletequation}\label{eqn:i2sum1}\sum_{i = 1}^{p}{\sum_{k = 2i - 1}^{2p}{i2^{-k}\binom{k+2}{2i+1}}} = p^2 + \frac{1}{2}p,  \end{bulletequation}

\item \begin{bulletequation}\label{eqn:i2sum2}\sum_{i = 1}^{p}{\sum_{k = 2i - 1}^{2p-1}{i2^{-k}\binom{k+2}{2i+1}}} = p^2 - \frac{1}{2}p \text{, and}\end{bulletequation}

\item \begin{bulletequation}\label{eqn:i2sum3}\sum_{i = 1}^{p+1}{\sum_{k = 2i - 1}^{2p+1}{i2^{-k}\binom{k+2}{2i+1}}} = p^2 + \frac{3}{2}p + \frac{1}{2}.\end{bulletequation}

\end{enumerate}
\end{lemma}

\begin{IEEEproof}
\begin{enumerate}[(i) ]

\item This follows by substituting $n = 2p + 4$ into \refe{eqn:binevensumi} and \refe{eqn:binevensum}, taking the difference between the first expression and twice the second, removing the $i = 0$ and $i = 1$ terms, shifting indices by $1$, and then dividing by $2$. Note that by the conditions on \refe{eqn:binevensumi} and \refe{eqn:binevensum}, this is true for $p \geq 0$.

\item This follows by substituting $n = 2p + 2$ into \refe{eqn:binoddsumi} and \refe{eqn:binoddsum}, taking the difference between these expressions and dividing by $2$. Note that by the conditions on \refe{eqn:binoddsumi} and \refe{eqn:binoddsum}, this is true for $p \geq 0$.

\item This follows by substituting $p+1$ for $p$ in \refe{eqn:ievenoddsum}.

\item This follows by substituting $p-1$ for $p$ in \refe{eqn:ievenevensum4}, adding this to \refe{eqn:ievenoddsum}, and using Pascal's rule to say that $\binom{2p+2}{2i+1} + \binom{2p+2}{2i+2} = \binom{2p+3}{2i+2}$.

\item This follows by substituting $n = 2p + 3$ into \refe{eqn:binoddsumi} and \refe{eqn:binoddsum}, taking the difference between these expressions, dividing by $2$ and taking the final term out of the sum. Note that by the conditions on \refe{eqn:binoddsumi} and \refe{eqn:binoddsum}, this is true for $p \geq 0$.

\item Let $s(p)$ represent the value of this sum, as a function of $p$. That is, \\ $s(p) = \sum_{i=1}^{p}{\sum_{k=2i-1}^{2p}{i2^{-k}\binom{k+2}{2i+1}}}$. Then we can see that $s(0) = 0$ and furthermore (using \refe{eqn:ioddoddsum} and \refe{eqn:ievenoddsum4}), $s(p+1) = s(p) +2(p+1) - \frac{1}{2}$. Thus, we can say that $s(p) = \sum_{k=1}^{p}{\left(2k-\frac{1}{2}\right)}$, which simplifies using \refe{eqn:sumint} to our desired result.

\item We can see that
\begin{multline*}
\sum_{i=1}^{p}{\sum_{k=2i-1}^{2p-1}{\!\!i2^{-k}\binom{k+2}{2i+1}}} = \sum_{i=1}^{p}{\sum_{k=2i-1}^{2p}{\!\!i2^{-k}\binom{k+2}{2i+1}}} \\
{} - 2^{-2p}\sum_{i=1}^{p}{i\binom{2p+2}{2i+1}},
\end{multline*}
and then the result follows from \refe{eqn:i2sum1} and \refe{eqn:ievenoddsum}.

\item This follows by substituting $p+1$ for $p$ in \refe{eqn:i2sum2}. \hfill \IEEEQEDhere

\end{enumerate}
\end{IEEEproof}

In addition to these series evaluations, the following series manipulations will prove to be useful.

\begin{lemma}[Equivalent binomial series]\label{lem:manips}
For integer values of $p \geq 0$ and $n \geq 0$,
\begin{enumerate}[(i) ]
\item  $\displaystyle\sum_{i = 1}^{p+1}{\!\sum_{k = 2i - 1}^{2p+1}{\!\!i2^{-k}\textstyle\binom{n+k+3}{n+2i+2}}} = \sum_{i = 1}^{p+1}{\!\sum_{k = 2i - 1}^{2p+1}{\!\!i2^{-k+1}\textstyle\binom{n+k+2}{n+2i+1}}} $
\begin{bulletequation}\label{eqn:i2manipodd}\hspace{3.5cm} {} - 2^{-2p-1}\sum_{i=1}^{p+1}{i\textstyle\binom{n+2p+4}{n+2i+2}},  \end{bulletequation} 
and

\item  $\displaystyle\sum_{i = 1}^{p}{\!\sum_{k = 2i - 1}^{2p}{\!\!i2^{-k}\textstyle\binom{n+k+3}{n+2i+2}}} = \sum_{i = 1}^{p}{\!\sum_{k = 2i - 1}^{2p}{\!\!i2^{-k+1}\textstyle\binom{n+k+2}{n+2i+1}}}$
\begin{bulletequation}\label{eqn:i2manipeven}\hspace{3.8cm}{} - 2^{-2p}\sum_{i=1}^{p}{i\textstyle\binom{n+2p+3}{n+2i+2}}.  \end{bulletequation}

\end{enumerate}
\end{lemma}

\begin{IEEEproof}
\begin{enumerate}[(i) ]
\item First, let us suppose that $p \geq 0$, $n \geq 0$ and $i$ is an integer between $1$ and $p+1$ (inclusive). Then, we can use Pascal's rule with $k > 2i - 1$ to write $\binom{n+k+3}{n+2i+2} = \binom{n+k+2}{n+2i+2} + \binom{n+k+2}{n+2i+1}$, while for $k = 2i - 1$ we can say $\binom{n+k+3}{n+2i+2} = \binom{n+2i+2}{n+2i+2} = 1 = \binom{n+2i+1}{n+2i+1} = \binom{n+k+2}{n+2i+1}$. With these two facts, we can write
\begin{multline}\label{eqn:proofmanipodd1}
\sum_{k = 2i-1}^{2p+1}{2^{-k}\binom{n+k+3}{n+2i+2}} = \sum_{k = 2i}^{2p+1}{2^{-k}\binom{n+k+2}{n+2i+2}} \\
	{} + \sum_{k = 2i-1}^{2p+1}{2^{-k}\binom{n+k+2}{n+2i+1}}.
\end{multline}
By shifting indices by $1$, substituting in \refe{eqn:proofmanipodd1}, and then rearranging, the first sum on the right becomes
\begin{multline*}
\sum_{k = 2i}^{2p+1}{2^{-k}\binom{n+k+2}{n+2i+2}} = \sum_{k = 2i-1}^{2p+1}{2^{-k}\binom{n+k+2}{n+2i+1}} \\
	{} - 2^{-2p-1}\binom{n+2p+4}{n+2i+2},
\end{multline*}
and so \refe{eqn:proofmanipodd1} becomes
\begin{multline*}
\sum_{k = 2i-1}^{2p+1}{2^{-k}\binom{n+k+3}{n+2i+2}} = 2\!\!\sum_{k = 2i-1}^{2p+1}{2^{-k}\binom{n+k+2}{n+2i+1}} \\
	{} - 2^{-2p-1}\binom{n+2p+4}{n+2i+2}.
\end{multline*}
Substituting this expression into the left hand side of \refe{eqn:i2manipodd} produces the desired result.

\item Again, let us suppose that $p \geq 0$, $n \geq 0$ and $i$ is an integer, now between $1$ and $p$ (inclusive). As above, we can use Pascal's rule to write
\begin{multline}\label{eqn:proofmanipeven1}
\sum_{k = 2i-1}^{2p}{2^{-k}\binom{n+k+3}{n+2i+2}} = \sum_{k = 2i}^{2p}{2^{-k}\binom{n+k+2}{n+2i+2}} \\
	{} + \sum_{k = 2i-1}^{2p}{2^{-k}\binom{n+k+2}{n+2i+1}}.
\end{multline}
By shifting indices by $1$, substituting in \refe{eqn:proofmanipeven1}, and then rearranging, the first sum on the right becomes
\begin{multline*}
\sum_{k = 2i}^{2p}{2^{-k}\binom{n+k+2}{n+2i+2}} = \sum_{k = 2i-1}^{2p}{2^{-k}\binom{n+k+2}{n+2i+1}} \\
	{} - 2^{-2p}\binom{n+2p+3}{n+2i+2},
\end{multline*}
and so \refe{eqn:proofmanipeven1} becomes
\begin{multline*}
\sum_{k = 2i-1}^{2p}{2^{-k}\binom{n+k+3}{n+2i+2}} = 2\!\!\sum_{k = 2i-1}^{2p}{2^{-k}\binom{n+k+2}{n+2i+1}} \\
	{} - 2^{-2p}\binom{n+2p+3}{n+2i+2}.
\end{multline*}
Substituting this expression into the left hand side of \refe{eqn:i2manipeven} produces the desired result. \hfill \IEEEQEDhere

\end{enumerate}
\end{IEEEproof}

Now, we can use Lemmas \ref{lem:specialsums} and \ref{lem:manips} to evaluate two more complicated expressions which will be necessary for the completion of our derivation.

\begin{lemma}\label{lem:oddexpression}
Let $p$ and $n$ be non-negative integers, and let
\begin{multline}\label{eqn:gdef}
g(n,p) \mathrel{\mathop :}= \frac{4p^2  \!+\! 6p \!+\! 2}{n\!+\!2p\!+\!2} \!+\! 4p \!+\! \left(4p^2 \!+\! 4np \!+\! 4n \!+\! 10p \!+\! 6\right)2^{1-n} \\
	{} + 2^{-2p}\\
	{} + 2^{-n-2p}\sum_{i=1}^{p+1} {\textstyle \!i\!\left\{2\binom{n+2p+4}{n+2i+1} \!-\! \binom{n+2p+4}{n+2i+2} \!-\! \left(2n\!+\!4p\!+\!6\right)\binom{2p+4}{2i+1}\right\}} \\
	{} + 2^{2-n}\sum_{i=1}^{p+1}{\sum_{k=2i-1}^{2p+1} {\textstyle\!\!\!\! i2^{-k}\!\left\{\frac{n+2p+1}{n+2p+2}\binom{n+k+2}{n+2i+1} \!-\! \binom{n+k+2}{n+2i} \!+\! \binom{k+3}{2i+1}\right\}}} \\
	{} + 2^{-2p}\sum_{i=1}^{p+1}{\sum_{k=1}^n {\textstyle\! i2^{-k}\left\{\binom{k+2p+4}{k+2i+2} \!-\! \frac{2n+4p+6}{n+2p+2}\binom{k+2p+3}{k+2i+1}\right\}}}.
\end{multline}
Then $g(n,p) = 0 \;\; \forall \; n,p \geq 0$.
\end{lemma}

\begin{IEEEproof}
First, we can use \refe{eqn:ievenoddsum4} to simplify the third term in the first sum. In addition, the third term in the second sum can be written as $\binom{k+2}{2i+1} + \binom{k+2}{2i}$ using Pascal's rule for $k \geq 2i-1$. We can then apply \refe{eqn:i2sum3} to the $\binom{k+2}{2i+1}$ term. This gives us
\begin{multline}\label{eqn:gsimp}
g(n,p) = \frac{4p^2\!  +\! 6p\! +\! 2}{n\!+\!2p\!+\!2} + 4p - \left(2p^2\! +\! 7p\! +\! 5\right)2^{1-n} + 2^{-2p} \\
	{} + 2^{-n-2p}\sum_{i=1}^{p+1} {i\textstyle\left\{2\binom{n+2p+4}{n+2i+1} - \binom{n+2p+4}{n+2i+2}\right\}} \\
	{} + 2^{2-n}\sum_{i=1}^{p+1}{\sum_{k=2i-1}^{2p+1} {\!\!\!\! i2^{-k}\left\{\textstyle\frac{n+2p+1}{n+2p+2}\binom{n+k+2}{n+2i+1} \!-\! \binom{n+k+2}{n+2i} \!+\! \binom{k+2}{2i}\right\}}} \\
	{} + 2^{-2p}\sum_{i=1}^{p+1}{\sum_{k=1}^n {\!i2^{-k}\left\{\textstyle\binom{k+2p+4}{k+2i+2} \!-\! \frac{2n+4p+6}{n+2p+2}\binom{k+2p+3}{k+2i+1}\right\}}}.
\end{multline}

Next, we will consider the case when $p = 0$. Using \refe{eqn:bincoef} and \refe{eqn:sumtwos}, we can simplify $g(n,0)$ to find that $g(n,0) = 0 \; \forall \; n \geq 0$. Thus, in the rest of the proof, we will assume that $p > 0$. Furthermore, when $n = 0$, we can use \refe{eqn:ievenevensum4}, \refe{eqn:ievenoddsum4} and \refe{eqn:i2sum3} to find that $g(0,p) = 0 \; \forall \; p > 0$.

Next, let us consider $g(n+1,p)$. Substituting $n+1$ in for $n$ in \refe{eqn:gsimp}, taking the $k = n+1$ terms out of the final sum and applying \refe{eqn:pascal} and \refe{eqn:i2manipodd} gives us
\begin{multline}\label{eqn:gn1p}
g(n\!+\!1,p)\! = \!\frac{4p^2\!  +\! 6p\! +\! 2}{n\!+\!2p\!+\!3} \!+\! 4p \!-\! \left(2p^2\! +\! 7p\! +\! 5\right)2^{-n} \!+\! 2^{-2p} \\
	{} + 2^{-n-2p-1}\sum_{i=1}^{p+1}{\! i\!\left\{\textstyle2\binom{n+2p+4}{n+2i+1} \!-\! 2\binom{n+2p+4}{n+2i+2}\right\}} \\
	{} + 2^{1-n}\sum_{i=1}^{p+1}{\sum_{k=2i-1}^{2p+1} {\!\!\!\! i2^{-k}\left\{\textstyle\frac{n+2p+1}{n+2p+3}\binom{n+k+2}{n+2i+1} \!-\! \binom{n+k+2}{n+2i} \!+\! \binom{k+2}{2i}\right\}}} \\
	{} + 2^{-2p}\sum_{i=1}^{p+1}{\sum_{k=1}^{n}{\! i2^{-k}\left\{\textstyle\binom{k+2p+4}{k+2i+2} \!-\! \frac{2n+4p+8}{n+2p+3}\binom{k+2p+3}{k+2i+1}\right\}}}.
\end{multline}

Now, let us define a new function, $a(n,p)$, as
\begin{equation}\label{eqn:adef}
a(n,p) \mathrel{\mathop :}= \left(n\!+\!2p\!+\!3\right)g(n\!+\!1,p) \!-\! \left(n\!+\!2p\!+\!2\right)g(n,p).
\end{equation}
Then, from \refe{eqn:gsimp} and \refe{eqn:gn1p}, we obtain
\begin{multline}\label{eqn:asimp}
a(n,p) = 4p \!+\! \left(4p^3 \!+\! 2np^2 \!+\! 16p^2 \!+\! 7np \!+\! 5n \!+\! 17p \!+\! 5\right)2^{-n} \\
	 	{} + 2^{-2p} + 2^{-n-2p}\sum_{i=1}^{p+1}{\textstyle\! i\left\{-\left(n\!+\!2p\!+\!1\right)\binom{n+2p+4}{n+2i+1} \!-\! \binom{n+2p+4}{n+2i+2}\!\right\}}  \\
		 {} + \left(n\!+\!2p\!+\!1\right)\!2^{1-n}\!\sum_{i=1}^{p+1}{\!\sum_{k=2i-1}^{2p+1} {\textstyle\!\!\!\! i2^{-k}\!\left\{\!-\binom{n+k+2}{n+2i+1} \!+\! \binom{n+k+2}{n+2i}\right. }}  \\
		 {\textstyle\left.{} \!-\! \binom{k+2}{2i}\!\right\}} \!+\! 2^{-2p}\sum_{i=1}^{p+1}{\sum_{k=1}^{n}{\textstyle\! i2^{-k}\!\left\{\!\binom{k+2p+4}{k+2i+2} \!-\! 2\binom{k+2p+3}{k+2i+1}\!\right\}}}.
\end{multline}
We can use \refe{eqn:ievenevensum4}, \refe{eqn:ievenoddsum4} and \refe{eqn:i2sum3} to show that $a(0,p) = 0$. In a similar manner as before, we will next consider $a(n+1,p)$. Substituting $n+1$ in for $n$ in \refe{eqn:asimp}, taking the $k = n+1$ terms out of the final sum  and applying \refe{eqn:pascal} and \refe{eqn:i2manipodd} produces
\begin{multline}\label{eqn:an1p}
a(n+1,p) \!=\! 4p +\! \left(4p^3 \!\!\!+\!\! 2np^2 \!\!\!+\! 18p^2 \!\!\!+\!\! 7np \!+\! 5n \!+\!\! 24p \!+\! 10\right)\!2^{-n-1} \\
	 	{} \!+\! 2^{-2p} \!+\! 2^{-n-2p-1}\!\sum_{i=1}^{p+1} {\textstyle\! i\!\left\{\!-\!\left(n\!+\!2p\!+\!2\right)\!\binom{n+2p+4}{n+2i+1} \!-\! 2\binom{n+2p+4}{n+2i+2}\!\right\}} \\
		{} + \left(n\!+\!2p\!+\!2\right)2^{-n}\!\sum_{i=1}^{p+1}{\!\sum_{k=2i-1}^{2p+1} {\textstyle\!\!\! i2^{-k}\left\{-\binom{n+k+2}{n+2i+1} + \binom{n+k+2}{n+2i}\right.}} \\
		{\textstyle\left.{} \!-\! \binom{k+2}{2i}\!\right\}} \!+\! 2^{-2p}\sum_{i=1}^{p+1}{\sum_{k=1}^{n}{\textstyle i2^{-k}\!\left\{\!\binom{k+2p+4}{k+2i+2} \!-\! 2\binom{k+2p+3}{k+2i+1}\!\right\}}}.
\end{multline}

Once again, we will define a new function, $b(n,p)$, as
\begin{equation}\label{eqn:bdef}
b(n,p) \mathrel{\mathop :}= \frac{a(n+1,p) - a(n,p)}{n+2p}.
\end{equation}
Note that $b(n,p)$ is well-defined since its denominator is positive for all $p > 0$ and $n \geq 0$. Then, from \refe{eqn:asimp} and \refe{eqn:an1p}, we obtain
\begin{multline}\label{eqn:bsimp}
b(n,p) = -\!\left(2p^2 \!+\! 7p \!+\! 5\right)2^{-n-1} + 2^{-n-2p-1}\!\sum_{i=1}^{p+1}{\textstyle i\binom{n+2p+4}{n+2i+1}}\\
	 {} \!-\! 2^{-n}\sum_{i=1}^{p+1}{\!\sum_{k=2i-1}^{2p+1} {\textstyle\!\!\!\! i2^{-k}\!\left\{\!-\binom{n+k+2}{n+2i+1} \!+\! \binom{n+k+2}{n+2i} \!-\! \binom{k+2}{2i}\!\right\}}}.
\end{multline}
Using \refe{eqn:ievenoddsum4} and \refe{eqn:i2sum3}, we find that $b(0,p) = 0$. Finally, we will follow our previous procedure once more and consider $b(n+1,p)$. Substituting $n+1$ in for $n$ in \refe{eqn:bsimp}, using \refe{eqn:pascal} and \refe{eqn:i2manipodd}, and comparing to \refe{eqn:bsimp} produces
\begin{equation}\label{eqn:bn1p}
b(n+1,p) = \frac{1}{2}b(n,p).
\end{equation}

Hence, from \refe{eqn:bn1p} and $b(0,p) = 0$, we conclude that $b(n,p) = 0 \;\; \forall \; n \geq0, \; p > 0$. Substituting this result into \refe{eqn:bdef} tells us that $a(n+1,p) = a(n,p) \;\; \forall \; n \geq0, \; p > 0$, which, along with the fact that $a(0,p) = 0$, allows us to conclude that $a(n,p) = 0 \;\; \forall \; n \geq0, \; p > 0$.

Finally, we can substitute this result into \refe{eqn:adef} to find that
\[
g(n+1,p) = \frac{n+2p+2}{n+2p+3}g(n,p)  \;\; \forall \; n \geq0, \; p > 0,
\]
which, along with the facts that $g(0,p) = 0$ and $g(n,0) = 0$, gives us our desired result.
\end{IEEEproof}

\begin{lemma}\label{lem:evenexpression}
Let $p$ and $n$ be non-negative integers, and let
\begin{multline}\label{eqn:hdef}
h(n,p) \mathrel{\mathop :}= \frac{4p^2  + 2p}{n\!+\!2p\!+\!1} + 4p - 2 + \left(4p^2 \!+\! 4np \!+\! 2n \!+\! 6p \!+\! 2\right)2^{1-n} \\
		{} + 2^{1-2p} - \left(4p^2 + 2np + 2n + 6p + 2\right)2^{1-n-2p}\\
		{} \!+\! 2^{1-n-2p}\!\sum_{i=1}^{p} {\textstyle\! i\!\left\{2\binom{n+2p+3}{n+2i+1} \!-\! \binom{n+2p+3}{n+2i+2} \!-\! \left(2n\!+\!4p\!+\!4\right)\!\binom{2p+3}{2i+1}\!\right\}} \\
		{} \!+\! 2^{2-n}\!\sum_{i=1}^{p}{\!\sum_{k=2i-1}^{2p} {\textstyle\!\!\!\! i2^{-k}\!\left\{\frac{n+2p}{n+2p+1}\binom{n+k+2}{n+2i+1} \!-\! \binom{n+k+2}{n+2i} \!+\! \binom{k+3}{2i+1}\!\right\}}} \\
		{} \!+\! 2^{1-2p}\!\sum_{i=1}^{p}{\!\sum_{k=1}^n {\textstyle\!\! i2^{-k}\!\left\{\binom{k+2p+3}{k+2i+2} \!-\! \frac{2n+4p+4}{n+2p+1}\binom{k+2p+2}{k+2i+1}\!\right\}}}.
\end{multline}
Then $h(n,p) = 0 \;\; \forall \; n,p \geq 0$.
\end{lemma}

\begin{IEEEproof}
This proof proceeds almost exactly as the proof of Lemma \ref{lem:oddexpression}. The only differences are that we use \refe{eqn:pascal}, \refe{eqn:ioddevensum}, \refe{eqn:ioddoddsum}, \refe{eqn:i2sum1} and \refe{eqn:i2manipeven} to simplify expressions (rather than \refe{eqn:sumtwos}, \refe{eqn:pascal}, \refe{eqn:ievenevensum4}, \refe{eqn:ievenoddsum4}, \refe{eqn:i2sum3} and \refe{eqn:i2manipodd}) and our intermediate functions are defined as
\[
c(n,p) \mathrel{\mathop :}= \left(n\!+\!2p\!+\!2\right)h(n+1,p) - \left(n\!+\!2p\!+\!1\right)h(n,p) \text{, and}
\]
\[
d(n,p) \mathrel{\mathop :}= \frac{c(n+1,p) - c(n,p)}{n+2p-1},
\]
where $d(n,p)$ is well-defined since its denominator is positive for all $p > 0$ and $n \geq 0$.
\end{IEEEproof}

Our final result covers some simplification required for the proof of Theorem \ref{theo:treeres}.
\begin{lemma}\label{lem:rsumsimp}
Suppose $n$ and $\ell$ are positive integers, and let
\begin{multline*}
s(n,\ell) \mathrel{\mathop:}= \frac{-3n^2 + 3\ell^2 - 2n\ell - n + 5\ell + 2}{2(n+\ell+1)^2} \\
		{} + \frac{\ell^2 + 2n\ell + 2n + 3\ell}{n+\ell+1}2^{-n} + \frac{n^2 + n + 2}{2(n+\ell+1)}2^{-\ell}\\
		{} + \tfrac{1}{2(n+\ell+1)} s_1(n,\ell) + \tfrac{1}{n+\ell+1} s_2(n,\ell) - \tfrac{n+\ell+2}{n+\ell+1} s_3(n,\ell) \\
		{} - \tfrac{n+\ell+2}{n+\ell+1} s_4(n,\ell) - \tfrac{1}{2(n+\ell+1)} s_5(n,\ell) - \tfrac{1}{n+\ell+1} s_6(n,\ell) \text{, where}
\end{multline*}
$s_1(n,\ell) \!\mathrel{\mathop :}=\! \displaystyle\sum_{k=1}^\ell {\!\left(4\! -\! \tfrac{2}{n+\ell+1}\! -\! 2^{k-\ell}\right)\!\left(n\!-\!k\right)}$, \newline
$s_2(n,\ell) \!\mathrel{\mathop :}=\! \displaystyle\sum_{k=1}^\ell {\!\left(4\! -\! \tfrac{2}{n+\ell+1}\! -\! 2^{k-\ell}\right)\!2^{1-n-k}\!\!\sum_{i=1}^{\left\lfloor\frac{k+1}{2}\right\rfloor}{\textstyle\!\! i\binom{n+k+2}{n+2i+1}}}$, \newline
$s_3(n,\ell) \mathrel{\mathop :}= \displaystyle\sum_{k=1}^n {\left(\tfrac{1}{n+\ell+1} - 2^{k-n}\right)\left(k-\ell\right)}$, \newline
$s_4(n,\ell) \mathrel{\mathop :}= \displaystyle\sum_{k=1}^n {\left(\tfrac{1}{n+\ell+1} - 2^{k-n}\right)2^{2-k-\ell}\sum_{i=1}^{\left\lfloor\frac{\ell+1}{2}\right\rfloor}{\textstyle i\binom{k+\ell+2}{k+2i+1}}}$, \newline
$s_5(n,\ell) \mathrel{\mathop :}= \displaystyle\sum_{k=1}^n {\sum_{j=1}^\ell {\left(2^{1+k-n} - 2^{j-\ell}\right)\left(k-j\right)}}$, and \newline 
$s_6(n,\ell) \mathrel{\mathop :}= \displaystyle\sum_{k=1}^n {\sum_{j=1}^\ell {\left(2^{1+k-n} - 2^{j-\ell}\right)2^{1-k-j}\sum_{i=1}^{\left\lfloor\frac{j+1}{2}\right\rfloor}{\textstyle i\binom{k+j+2}{k+2i+1}}}}$.

Then
\begin{multline}\label{eqn:rsumsimp}
s(n,\ell) = 2\left(n - \ell - 1\right) + 2^{2-n-\ell}\sum_{i=1}^{\left\lfloor\frac{\ell+1}{2}\right\rfloor}{ i\binom{n+\ell+3}{n+2i+1}} \\
			{} + \frac{1}{n+\ell+1}\left[\frac{\ell^2 + \ell}{n+\ell+1} + 2\ell - 2 + 2^{1-\ell} \vphantom{\sum_{i=1}^{\left\lfloor\frac{\ell+1}{2}\right\rfloor}{i}}\right. \\
			{} + \left(\ell^2 + 2n\ell + 2n + 3\ell + 2\right)2^{1-n} \\
			{} + 2^{1-\ell}\sum_{i=1}^{\left\lfloor\frac{\ell+1}{2}\right\rfloor}{\sum_{k=1}^{n}{\textstyle i2^{-k}\left\{\binom{k+\ell+3}{k+2i+2} - \frac{2n+2\ell+4}{n+\ell+1}\binom{k+\ell+2}{k+2i+1}\right\}}} \\
			{} \!+\! 2^{1-n-\ell}\!\!\sum_{i=1}^{\left\lfloor\frac{\ell+1}{2}\right\rfloor}{\textstyle\!\!\! i\!\left\{2\binom{n+\ell+3}{n+2i+1} \!-\! \binom{n+\ell+3}{n+2i+2} \!-\! \left(2n \!+\! 2\ell\!+\!4\right)\binom{\ell+3}{2i+1}\!\right\}} \\
			\left.{} \!+\! 2^{2-n}\!\!\!\!\sum_{i=1}^{\left\lfloor\frac{\ell+1}{2}\right\rfloor}{\!\!\!\!\sum_{k=2i-1}^{\ell}{\textstyle\!\!\!\! i2^{-k}\!\left\{\!\frac{n+\ell}{n+\ell+1}\binom{n+k+2}{n+2i+1} \!-\! \binom{n+k+2}{n+2i} \!+\! \binom{k+3}{2i+1}\!\right\}}}\! \right].
\end{multline}
\end{lemma}

\begin{IEEEproof}
We can use \refe{eqn:sumint}, \refe{eqn:sumtwos} and \refe{eqn:suminttwos} to simplify each of $s_1(n,\ell)$, $s_3(n,\ell)$ and $s_5(n,\ell)$. This gives us $s_1(n,\ell) = \tfrac{-2\ell^3 + 4n^2\ell+2n\ell^2-2n^2-\ell^2-4n-\ell-2}{n+\ell+1} + \left(n+1\right)2^{1-\ell}$, $s_3(n,\ell) = \tfrac{-3n^2+4\ell^2-2n\ell+n+8\ell+4}{2\left(n+\ell+1\right)} - \left(\ell+1\right)2^{1-n}$ and $s_5(n,\ell) = -n^2 - 2\ell^2 + 6n\ell - 3n - 6\ell + \left(\ell^2 + 3\ell\right)2^{1-n} + \left(n^2 + 3n\right)2^{-\ell}$.

To simplify each of $s_2(n,\ell)$, $s_4(n,\ell)$ and $s_6(n,\ell)$, we can first exchange the order of summation to make the sum over $i$ the outermost sum, and then apply \refe{eqn:pascal}, \refe{eqn:binsumnm} and/or \refe{eqn:binsumnmk} to obtain
\begin{multline*}
s_2(n,\ell) = \left(4 - \tfrac{2}{n+\ell+1}\right)2^{1-n}\sum_{i = 1}^{\left\lfloor\frac{\ell+1}{2}\right\rfloor}{\sum_{k=2i-1}^{\ell}{\textstyle i2^{-k}\binom{n+k+2}{n+2i+1}}}\\
{} - 2^{1-n-\ell}\sum_{i = 1}^{\left\lfloor\frac{\ell+1}{2}\right\rfloor}{\textstyle i\binom{n+\ell+3}{n+2i+2}},
\end{multline*}
\begin{multline*}
s_4(n,\ell) = \tfrac{1}{n+\ell+1}2^{2-\ell}\sum_{i = 1}^{\left\lfloor\frac{\ell+1}{2}\right\rfloor}{\sum_{k=1}^{n}{\textstyle i2^{-k}\binom{k+\ell+2}{k+2i+1}}}\\
{} - 2^{2-n-\ell}\sum_{i = 1}^{\left\lfloor\frac{\ell+1}{2}\right\rfloor}{\textstyle i\binom{n+\ell+3}{n+2i+1}} + 2^{2-n-\ell}\sum_{i = 1}^{\left\lfloor\frac{\ell+1}{2}\right\rfloor}{\textstyle i\binom{\ell+3}{2i+1}} \text{, and}
\end{multline*}
\begin{multline*}
s_6(n,\ell) = 2^{2-n}\sum_{i=1}^{\left\lfloor\frac{\ell+1}{2}\right\rfloor}{\!\sum_{k=2i-1}^{\ell}{\textstyle\! i2^{-k}\!\left\{\binom{n+k+2}{n+2i+1} + \binom{n+k+2}{n+2i}\right\}}} \\
{} - 2^{2-n}\!\sum_{i=1}^{\left\lfloor\frac{\ell+1}{2}\right\rfloor}{\!\!\sum_{k=2i-1}^{\ell}{\textstyle\!\! i2^{-k}\binom{k+3}{2i+1}}} - 2^{1-\ell}\!\sum_{i=1}^{\left\lfloor\frac{\ell+1}{2}\right\rfloor}{\!\sum_{k=1}^{n}{\textstyle i2^{-k}\binom{k+\ell+3}{k+2i+2}}}.
\end{multline*}

Substituting our simplified expressions for $s_1, s_2, s_3, s_4, s_5$ and $s_6$ into the definition of $s(n,\ell)$ gives us our desired result.
\end{IEEEproof}


\bibliographystyle{IEEEtran}
\bibliography{REFabrv,ReferenceList}

\end{document}